\documentclass[a4paper,10pt]{article}
\usepackage[top=30truemm,bottom=25truemm,left=28truemm,right=22truemm]{geometry}
\usepackage[dvipdfmx,hiresbb]{graphicx}
\usepackage{wrapfig}
\usepackage{here}
\usepackage{framed}
\usepackage{amsmath,ascmac,amsfonts,amssymb,amsthm,pgfpages,enumerate,cases}
\usepackage{algorithm}
\usepackage{algorithmic}
\usepackage{color}
\usepackage{titlesec}
\usepackage{listings}
\usepackage[english]{babel}
\definecolor{lightgray}{rgb}{0.75,0.75,0.75}
\definecolor{darkgray}{rgb}{0.25,0.25,0.25}
\usepackage{bm}

\newcommand{\R}{\mathbb{R}}

\newcommand{\x}{\bm{x}}
\newcommand{\y}{\bm{y}}

\renewcommand{\u}{\bm{u}}
\renewcommand{\v}{\bm{v}}
\newcommand{\f}{\bm{f}}
\newcommand{\n}{\bm{n}}

\renewcommand{\div}{\mathrm{div}}

\title{A variational approach to the \\inverse imaging of composite elastic materials}
\author{Elliott Ginder\\ {\emph{School of Interdisciplinary Mathematical Sciences}}\\{\emph{Meiji University}} \and Riku Kanai\\ {\emph{Graduate School of Advanced Mathematical Sciences}}\\{\emph{Meiji University}}}
\date{}
\begin{document}
\maketitle
\begin{abstract}
We introduce a framework for performing the inverse imaging of composite elastic materials. Our technique uses surface acoustic wave (SAW) boundary observations within a minimization problem to express the interior composition of the composite elastic materials. We have approached our target problem by developing mathematical and computational methods for investigating the numerical solution of the corresponding inverse problem. 

We also discuss a mathematical model for expressing the propagation of elastic
waves through composite elastic bodies, and develop approximation schemes for investigating its numerical solutions.
Using these methods, we define a cost functional for measuring the difference between simulated and given SAW data. Then, using a Lagrangian approach, we are able to determine the gradient of the cost functional and analyze the inverse imaging problem's solution as a gradient flow. The cost functional's gradient is composed of solutions to a state equation, as well as of solutions to related adjoint problems. We thus developed numerical methods for solving these problems and investigated the gradient flow of the cost functional. Our results show that the gradient flow is able to recover the interior composition of the composite, and we illustrate this fact using the numerical realization of our proposed framework.
\end{abstract}
\maketitle

\section{Introduction}
Designing accurate simulations of natural and social phenomena undoubtably requires the effective integration of quantitative and qualitative information. Such information is often gained through both theoretic and experimental methods, and techniques for its incorporation represent an important theme in data assimilation. It is also well known that relating theoretical assumptions and experimental observations to simulations often leads to the mathematical formulation of an inverse problem (see e.g., \cite{Nakamura}). From CT (computed tomography) scanning techniques in the medical sciences, to sonic cloaking and sonar technology, inverse problems abound in applied mathematics and engineering. 

In \cite{Ohtsuka}, the authors designed an experiment for obtaining {\emph{surface acoustic wave} (SAW)} data of a phononic crystal and used k-space imaging techniques to analyze their experimental results. In addition to estimating band gaps of the phononic crystal, the authors showed that the crystal's approximate surface geometry can be imaged using the acoustic field data. The current research investigates the possibility of using SAW data to image not only the exterior of composite media, but the interior as well. In particular, we design and analyze an inverse problem for imaging the interior of composite elastic materials, using SAW data.

%
%

This paper is outlined as follows. A mathematical model for describing elastic wave propagation through composite media is introduced in section \ref{model}. This section also discusses a level-set approach for describing the location of the different materials. Section \ref{approx} briefly describes the numerical approximation of solutions to the model equation, and these techniques are used in section \ref{numbehave} to illustrate the target phenomena's simulation. The main target of this research, which is an inverse problem for describing the interior of composite elastic materials, is introduced in section \ref{inv}. The inverse problem is formulated as a minimization between boundary (experimental) data and simulated data. We determine the gradient of the corresponding cost functional in section \ref{lag}, and a level set method for approximating solutions of the inverse problem is formulated in section \ref{level}. The numerical results are shown in section \ref{numresult}, where we show that the gradient flow of the cost function indeed contains information on the interior of the composite media. We give our conclusions, including areas for future work, in section \ref{conclusion}.
\vspace{70pt}

\section{Modeling of wave propagation through composite elastic materials}\label{model}

This section describes the mathematical model equation used in simulating the propagation of elastic waves through composite elastic media. We also define a level set framework for expressing the location of material parameters. We will discuss our method for constructing numerical solutions to the model equations and show the numerical behavior of the model equation under an induced outer force. 

\begin{figure}[H]
  \centering
  \includegraphics[scale=0.5]{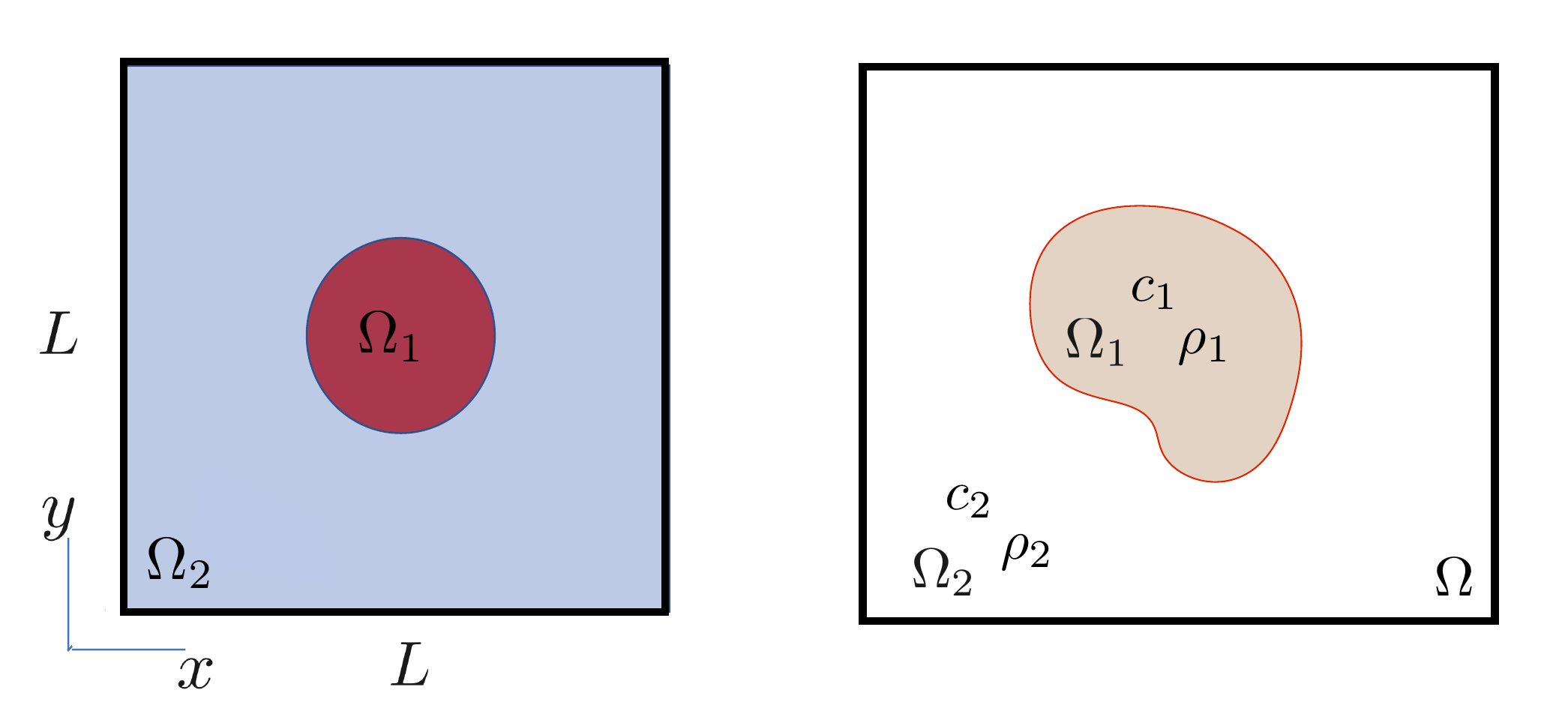}
  \caption{Left: a two dimensional representation of a composite elastic material. Right: spatial dependence of elastic parameters.}
  \label{fig:modeling1}
\end{figure}

We represent the composite media as a bounded, smooth, open subset $\Omega$ of ${\bf{R}^3}$ or ${\bf{R}^2}$ which is formed by a union:
\begin{align}
\Omega = \Omega_1\cup \Omega_2.
\end{align}
Here $\Omega_1$ and $\Omega_2$ are themselves smooth, open subsets of the ambient space which represent the two media (see figure \ref{fig:modeling1} for a schematic representation in the two dimensional setting).  In what follows, for the sake of clarity, we will assume the three dimensional setting, writing $\x = (x,y,z)$.

\begin{figure}[H]
  \centering
  \includegraphics[scale=0.25]{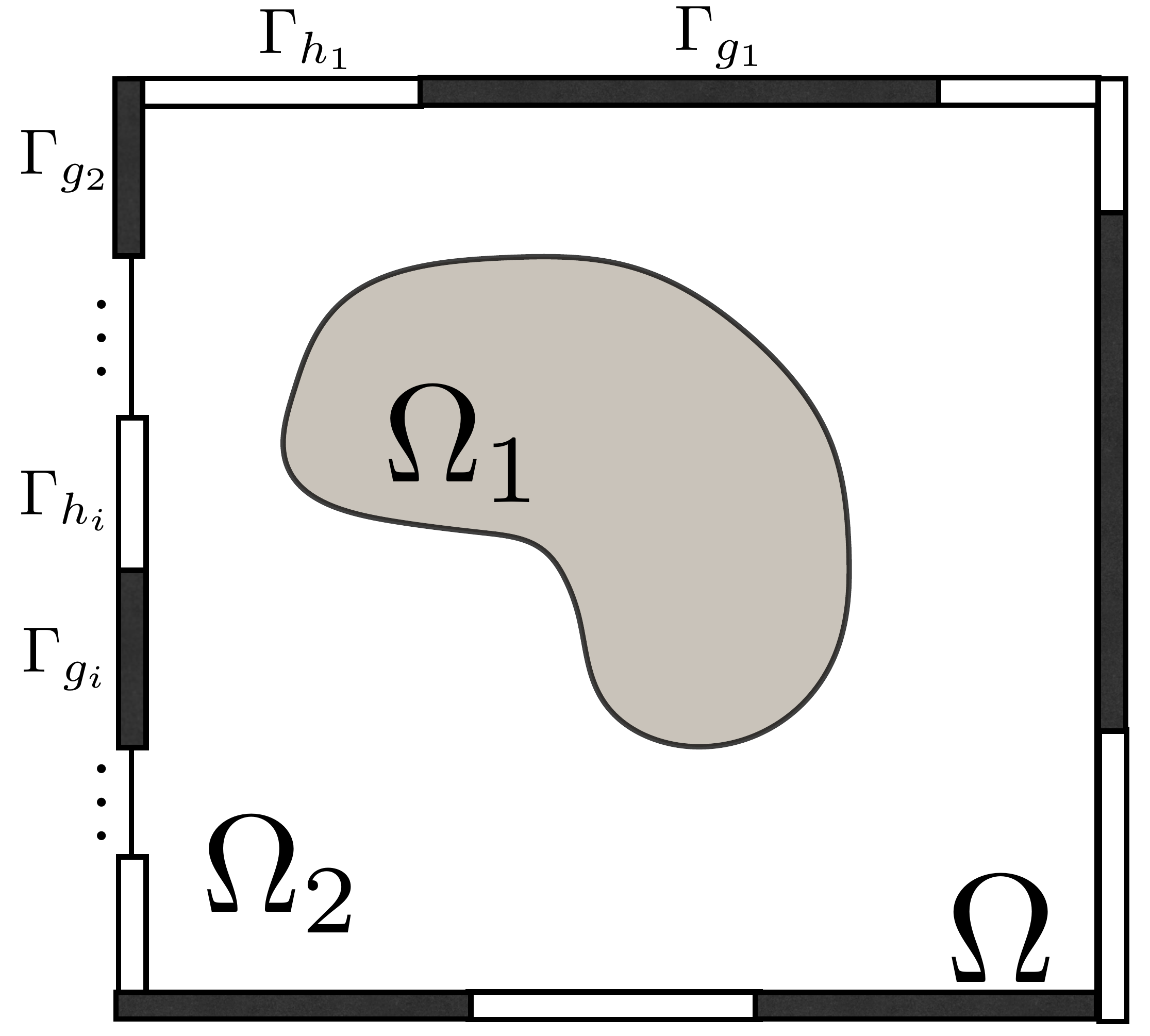}
  \caption{Schematic representation of the exterior boundary conditions.}
  \label{fig:bc}
\end{figure}

The elastic wave propagation is then assumed to be described by the displacement vector field $\u:[0,T)\times\Omega\rightarrow {\bf{R}^3}$, which evolves as a solution to the composite elastic wave equation (see e.g., \cite{Banerjee} or \cite{Laude}):
\begin{eqnarray}
  \begin{cases}
    \rho(\x)\u_{tt} = \div(c(\x)\otimes\nabla\u) + \f(\x,t)
    & \text{in}\;(0,T)\times\Omega\\
    u_i = g_i & \text{on}\;(0,T)\times\Gamma_{g_i}\;(\text{Dirichlet B.C.})\\
    (c(\x)\otimes\nabla\u)_{ij}n_j = h_i
    & \text{on}\;(0,T)\times\Gamma_{h_i}\;(\text{Neumann B.C.})\\
    [u_i(\x,t)] = 0 & \text{on}\;(0,T)\times\partial\Omega_1\cap\partial\Omega_2\\
    [(c(\x)\otimes\nabla\u)_{ij}] = 0 & \text{on}\;(0,T)\times\partial\Omega_1\cap\partial\Omega_2\\
    u_i(0,\x) = u_{0i}(\x) & \x \in \Omega\\
    u_{i,t}(0,\x) = \dot{u}_{0i}(\x) & \x\in\Omega.
  \end{cases}
  \label{eq:composite wave}
\end{eqnarray}
When $\Omega$ is a subset of ${\bf{R}^3}$ the indices take the values $i=1,2,3$ and $j=1,2,3$. Similarly, whenever $\Omega$ is a subset of ${\bf{R}^2}$, the indices take the values $i=1,2$ and $j=1,2$. We also formally assume 
that the exterior of $\Omega$:
\begin{align}
	\partial\Omega \backslash \{\partial \Omega_1 \cap \partial\Omega_2\} = \left\{{\bigcup_{i}} \Gamma_{h_i} \right\} \cup \left\{{\bigcup_{i}} \Gamma_{g_{i}}\right\}
\end{align}
is composed of regions $\Gamma_{g_i}$ and $\Gamma_{h_i}$, where Dirichlet and Neumann boundary conditions are prescribed, repectively (see figure, \ref{fig:bc}). In addition to the exterior boundary conditions on $\Omega$, interfacial boundary conditions holding on boundaries between $\Omega_1$ and $\Omega_2$ are prescribed. These are expressed as jump conditions, where $[\cdot]$ denotes the jump of its quantity across the normal to the interface (see e.g., \cite{Laude})．
In the above, $\nabla\u$ refers to the matrix:
\begin{subequations}
  \begin{align}
    \nabla\u =
    \begin{pmatrix}
      \dfrac{\partial u_1}{\partial x} & \dfrac{\partial u_1 }{\partial y} &  \dfrac{\partial u_1 }{\partial z}\\
      \dfrac{\partial u_2}{\partial x} & \dfrac{\partial u_2 }{\partial y} &  \dfrac{\partial u_2 }{\partial z}\\
      \dfrac{\partial u_3}{\partial x} & \dfrac{\partial u_3 }{\partial y} &  \dfrac{\partial u_3 }{\partial z}\notag
    \end{pmatrix},
  \end{align}
\end{subequations}
$c(\x)$ is a spatially dependent stiffness tensor, and $\rho(\x)$ is a spatially dependent mass density function. The notation $c(\x)\otimes\nabla\u$ denotes:
\begin{equation}
  c(\x)\otimes\nabla\u = \sum_{i,j=1}^3 c(\x)_{ijkl}\nabla\u_{ij}.
\end{equation}

\begin{figure}[H]
  \centering
  \includegraphics[scale=0.5]{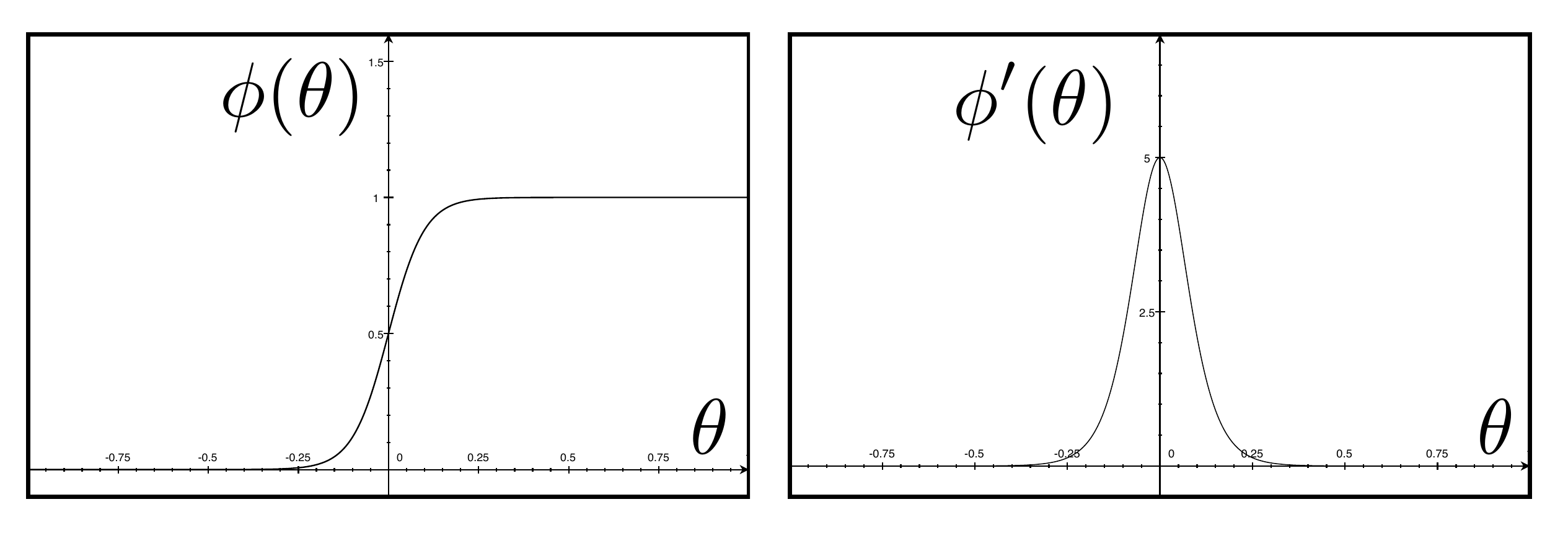}
  \caption{Left: the sigmoid function $\phi(\theta)$. Right: its derivative, ($\varepsilon = 0.1$).}
  \label{fig:sig}
\end{figure}

The spatial dependence of $\rho(\x)$ and $c(\x)$ are described using a signed distance level set function:
\begin{align}
  \theta(\x) =
  \begin{cases}
    -\text{dist}_{\partial\Omega_1}(\x) & \x \in \Omega_1\\
    \text{dist}_{\partial\Omega_1}(\x) & \x \in \Omega_2,\\
  \end{cases}
  \label{eq:levelset}
\end{align}
where $\text{dist}_{E}(\x)$ is defined to be the distance from the point $\x$ to the set $E$:
\begin{align}
 \text{dist}_{E}(\x) =
  \inf_{y\in E}||\x-\y|| \hspace{20pt} \x \in \Omega.
  \label{eq:dist}
\end{align}
Here $||\cdot||$ denotes the euclidean norm. We remark that the level set function $\theta(\x)$ allows us to determine the location of $\Omega_1$ and $\Omega_2$ as follows:
\begin{eqnarray}
  \Omega_{1} := \{\x\in\Omega|\theta(\x)<0\}\\
  \Omega_{2} := \{\x\in\Omega|\theta(\x)\ge 0\}.
\end{eqnarray}
The level set function is given as the argument of a sigmoid function (see figure \ref{fig:sig}): 
\begin{equation}
  \phi(\theta) = \dfrac{1}{2}\left(\tanh\left(\dfrac{\theta}{\varepsilon}\right)+1\right),
\end{equation}
where $0< \varepsilon\ll 1$ is an interpolation parameter.
\noindent This allows us to control the mass density and elastic tensors using the interpolation function:
\begin{eqnarray}\label{eq:parameters}
  \rho(\theta) &=& \rho_1\phi(\theta) + \rho_2(1 - \phi(\theta)),\notag\\
  c(\theta) &=& c_1\phi(\theta) + c_2(1 - \phi(\theta)),
\end{eqnarray}
where $\rho_1,\rho_2$ and $c_1,c_2$ correspond to the parameters within $\Omega_1$ and $\Omega_2$, respectively.
\noindent For later use in our numerical computations, we remark that the derivative of the sigmoid function is an approximation of a
delta function (see figure \ref{fig:sig}):
\begin{equation}
  \dfrac{d \phi(\theta)}{d \theta} = \dfrac{1}{2\varepsilon}\text{sech}^2\left(\dfrac{\theta}{\varepsilon}\right).
\end{equation}
The term $\f(\x,t)$ in equation (\ref{eq:composite wave}) denotes a time dependent outer force, whereas $u_{0i}(\x)$ and $\dot{u}_{0i}(\x)$ designate the initial displacement and velocity fields.


\section{Approximation of solutions to the model equation} \label{approx}
This section briefly discusses the methods used in constructing numerical solutions of the composite elastic wave equation ($\ref{eq:composite wave}$). 


Time is discretized using a time step $h>0$ and we have:
\begin{equation}
  \dfrac{\partial^2 \u}{\partial t^2} = \dfrac{\u - 2\u_{n} + \u_{n-1}}{h^2} +O(h^2),
  \hspace{15pt}\text{(as $h\rightarrow 0$)}
\end{equation}
where $\u_n$ and $\u_{n-1}$ are given vector valued functions approximating $\u(\x,nh)$ and $\u(\x,(n-1)h)$, respectively. That is, for suitably small $h>0,$ we have the approximation:
\begin{equation}\label{timeapprox}
  \dfrac{\partial^2 \u}{\partial t^2} \approx \dfrac{\u - 2\u_{n} + \u_{n-1}}{h^2},
\end{equation}
where $n$ is a nonnegative integer. When $n=0$, $\u_0(\x)$ denotes the initial displacement field and $\u_{-1}$ is constructed using the initial velocity field $\v_0(\x)$ within a backward differencing:
\begin{align}
  \u_{-1} = \u_0 - h\v_0.
\end{align}
Inserting approximation (\ref{timeapprox}) into the model equation (\ref{eq:composite wave}) and rearranging yields the following elliptic partial differential equation for $\u(\x)$:
\begin{equation}
  \rho(\theta)\u = h^2\nabla\cdot(c(\theta) \otimes \nabla\u) + h^2\f_n + \rho(\theta)
  (2\u_{n}-\u_{n-1}).
  \label{eq:sabun}
\end{equation}
Here $\f_n(x)$ denotes $\f(\x,nh)$. For prescribed boundary conditions, we then define $\u_n$ to be the solution of equation ($\ref{eq:sabun}$) and are thus able to approximate the solution of the model equation ($\ref{eq:composite wave}$) as the sequence of functions $\{\u_n\}_{n=0}^{N},$ where $N$ denotes the final time step.

The numerical realization of each $\u_n$ is obtained as a weak solution of ($\ref{eq:sabun}$). In particular, assuming Dirichlet zero boundary conditions for simplicity, we multiply equation $(\ref{eq:sabun})$ by an arbitrary test function $\v\in C^{\infty}_0(\Omega;\R^d)$ and integrate over $\Omega$:
\begin{eqnarray}
  \int_\Omega \rho(\theta)\u\cdot\v d\x &=&
  \int_\Omega ((h^2\nabla\cdot(c(\theta)\otimes\nabla\u)) + h^2\f_n + \rho(\theta)(2\u_{n}-\u_{n-1}) )\cdot\v d\x,\\
  &=& h^2 \int_\Omega \nabla\cdot(c(\theta)\otimes\nabla\u)\cdot\v d\x + \int_\Omega (h^2\f_n + \rho(\theta)(2\u_{n}-\u_{n-1}))\cdot\v d\x.
\end{eqnarray}

Applying Gauss-Green's theorem (see e.g., \cite{Pironneau}) and using the boundary conditions expresses
\begin{equation*}
  h^2\int_\Omega \nabla\cdot(c(\theta)\otimes\nabla\u)\cdot\v d\x =
  -h^2\int_\Omega(c(\theta)\otimes\nabla\u):\nabla\v d\x.
  \label{eq1}
\end{equation*}
The weak form of (\ref{eq:sabun}) is then expressed:
\begin{equation}
  \int_\Omega \rho(\theta)\u\cdot\v d\x =
  -h^2\int_\Omega(c(\theta)\otimes\nabla\u):\nabla\v + h^2f\cdot\v + \rho(\theta)(2\u_{n}-\u_{n-1})\cdot\v d\x,
  \label{weak1}
\end{equation}
for all test functions $\v$ in $C^{\infty}_{0}(\Omega; {\bf{R}}^d)$. 

Using the Lagrange $P^1$ finite element method (see e.g., \cite{Hughes}), we are thus able to obtain
each $\{\u_{n}\}_{n=0}^N$ as the solution of a system of linear equations stemming from ($\ref{weak1}$).

\section{Numerical behavior of the model equation}\label{numbehave}
\begin{figure}[H]
  \begin{minipage}{0.5\hsize}
    \begin{center}
      \includegraphics[scale=0.45]{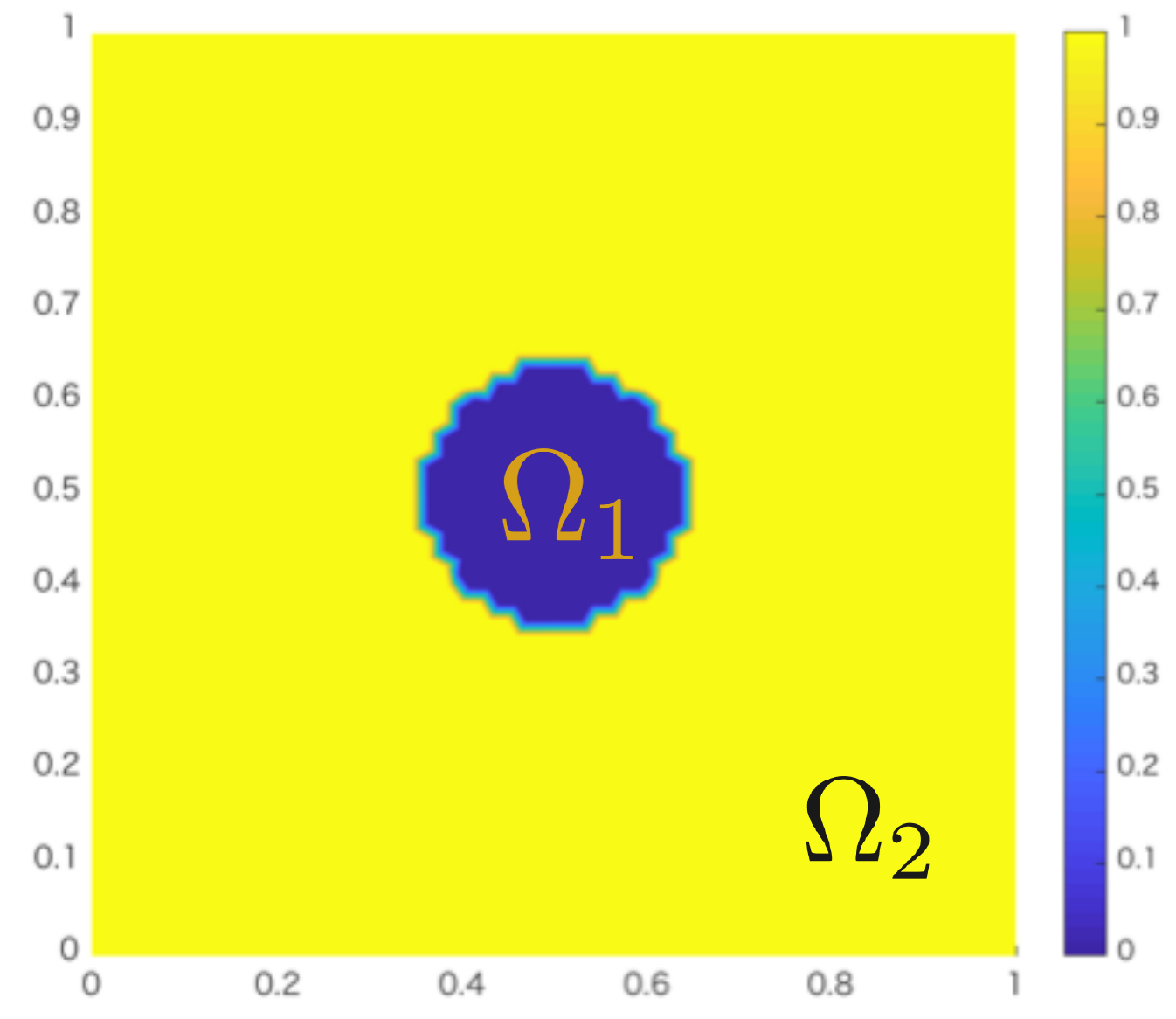}
    \end{center}
    \caption{Truncation of the level set function, $\chi_{\Omega_2}(\x)$.}
    \label{theta}
  \end{minipage}
  \begin{minipage}{0.5\hsize}
    \begin{center}
      \includegraphics[scale=0.4,trim=0cm 0.2cm 0cm 0cm,clip=true]{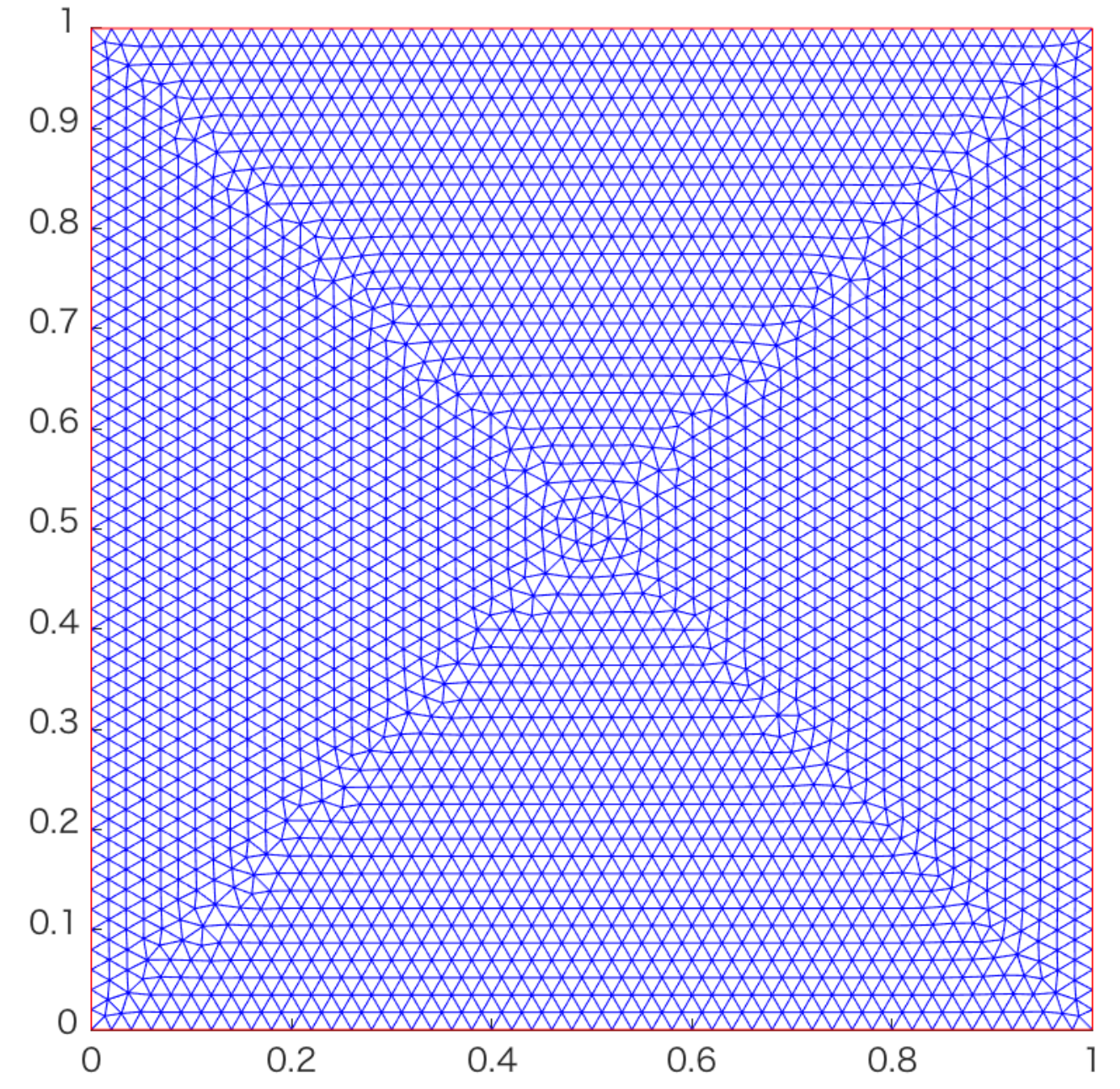}
    \end{center}
    \caption{The finite element mesh.}
    \label{Delaunay}
  \end{minipage}
\end{figure}

Using the approximation method described in the previous section, we will now briefly illustrate the numerical behavior of solutions to the model equation (\ref{eq:composite wave}) in the two dimensional setting. A Delaunay triangulation is employed on 5266 nodes to partition the domain $\Omega=(0,1)\times(0,1)$ into a finite number of elements, each with a target maximum edge length of 0.015 (figure \ref{Delaunay} illustrates the mesh). We prescribe a time step approximately equal to $h=3\times10^-6$.

The mass density function and stiffness tensors are prescribed using equation ($\ref{eq:parameters}$), where $\varepsilon=0.05$ and 
\begin{align}
\partial \Omega_1=\{\x\in {\bf{R}}^2: ||\x-(0.5,0.5) || = 0.1\}.
\label{eq:location}
\end{align}
The corresponding level set function is given by equation ($\ref{eq:levelset}$) and figure \ref{theta} shows the truncation of the level set function:
\begin{align}
	\chi_{\Omega_2}(\x)=\begin{cases}
	1\hspace{20pt} \x \in \Omega_2\\
	0\hspace{20pt} \text{otherwise.}
	\end{cases}\notag
\end{align}

\begin{figure}[H]
  \centering
  \includegraphics[scale=0.5]{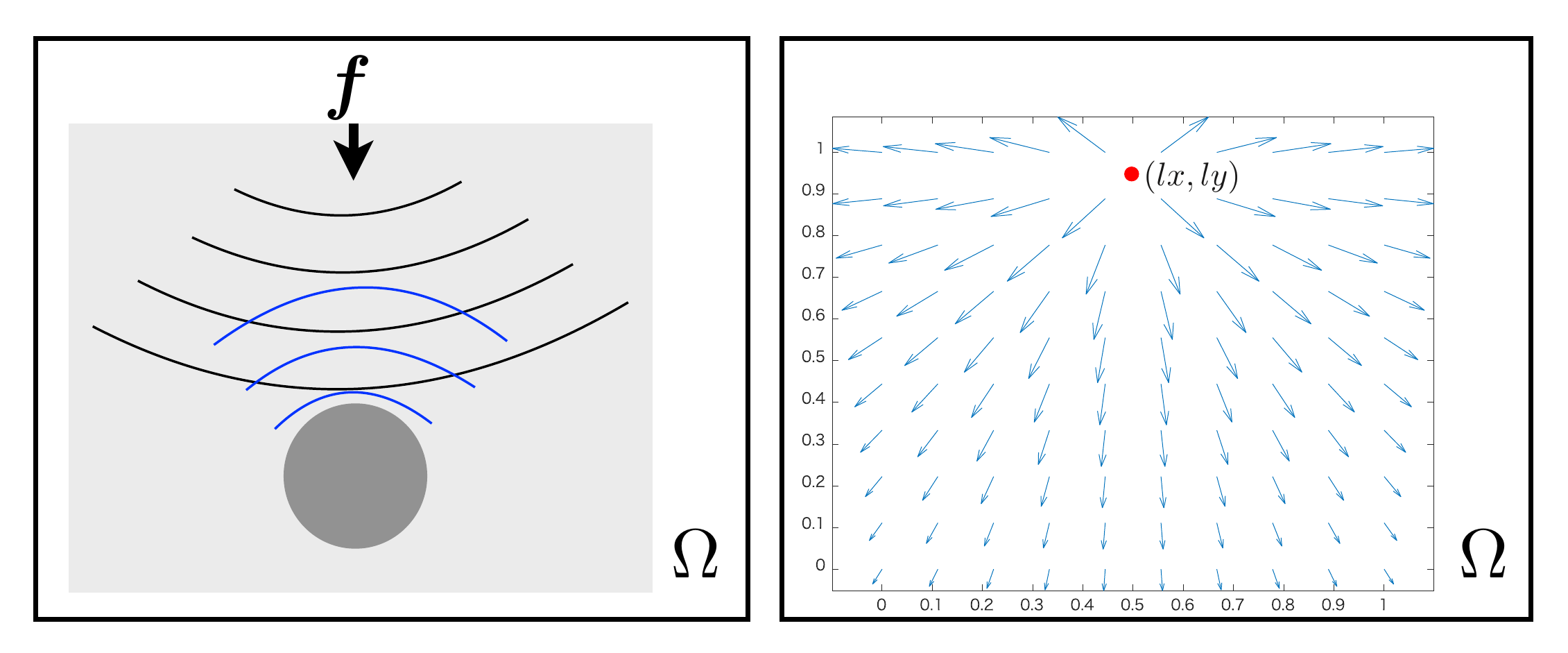}
  \caption{Left: a schematic representation of the boundary outer forcing. Right: a representative outer forcing in equation (\ref{mylaser}), for certain parameter values.}
  \label{fig:outerforce}
\end{figure}

The purpose of the following computations is to show the influence of the interior structure on the wave propagation. For our first computation, the domain is assumed to be composed of two isotropic materials. For the sake of illustration, we prescribe $\Omega_2$ to represent the region of a glass and $\Omega_1$ to represent a steel region. Omitting units, the stiffness tensors are then expressed in terms of their Poisson's ratio and Young's modulus:
\begin{align}
&\nu_1 = 0.26,\hspace{30pt} E_1 = 180\times 10^9\\
&\nu_2 = 0.25,\hspace{30pt} E_2 = 70\times 10^9.
\end{align}
Correspondingly, the density functions are set with the values:
\begin{align}
&\rho_1 = 4\times10^3\\
&\rho_2 = 8\times10^3.
\end{align}
Elastic waves are generated by imparting an outer force $\f(\x,t)$ near the top boundary of the domain, described as follows:
\begin{table}[H]
  \centering
  \begin{tabular}{ll}
    $\f(\x,t)=(f_1,f_2)^T$ &: a time dependent outer force in the $(x_1,x_2)^T$ direction\\
    $(l_x,l_y)=(0.5,0.98)$ &:  central focus of the outer force\\
    $l_N=4$& :  number of oscillations within the time interval $(0,T)$\\
    $A=10^{10}$& : a parameter controlling the magnitude of the outer force.\\
    $l_{width}=1300$&  : a parameter controlling the extent of the outer force.\\
  \end{tabular}
\end{table}
\noindent We also set
\begin{align}
  r(\x)= \sqrt{(x-l_x)^2 + (y-l_y)^2)} \hspace{15pt} \text{ and } \hspace{15pt} 
  \psi(\x) = \arctan\left(\dfrac{y-l_y}{x-l_x}\right).
\end{align}
The outer force used in our simulation is then defined (the vector field in figure $\ref{fig:outerforce}$ represents an instant in time of one such outer force):
\begin{equation}
  \f(\x,t) = A\cos\left(\frac{2\pi t}{T}l_{N}\right)\exp(-l_{width} r)
  \begin{pmatrix}
    \cos\psi \\
    \sin\psi
  \end{pmatrix}.
  \label{mylaser}
\end{equation}

\begin{figure}[H]
  \centering
  \includegraphics[page=1,scale=0.6]{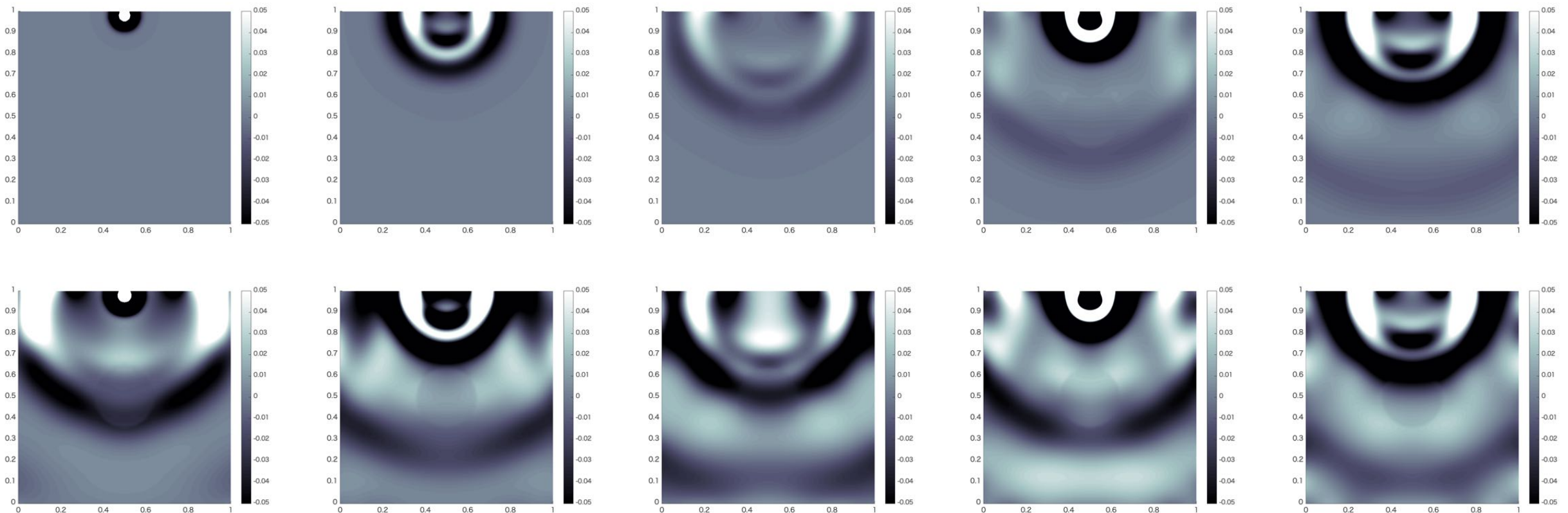}
  \caption{Evolution of the dilation of $\u$, with the materials set as in figure (\ref{theta}). Time is from top to bottom, left to right.}
  \label{u}
\end{figure}

Using the approximation method described in the previous section, we solve the model equation ($\ref{eq:composite wave}$) for $\u$ under Neumann boundary conditions, until approximately the time $T=4\times10^{-4}$. The numerical results are shown in figure \ref{u}, where we display the dilation of the displacement field:
\begin{equation}
  \text{dilation}\;\u = \dfrac{\div\u}{|\nabla\u|}.
\end{equation}
Close inspection of figure \ref{u} shows the spatial dependence imparted by the inclusion $\Omega_1$. 

\begin{figure}[H]
  \centering
  \includegraphics[page=2,scale=0.6]{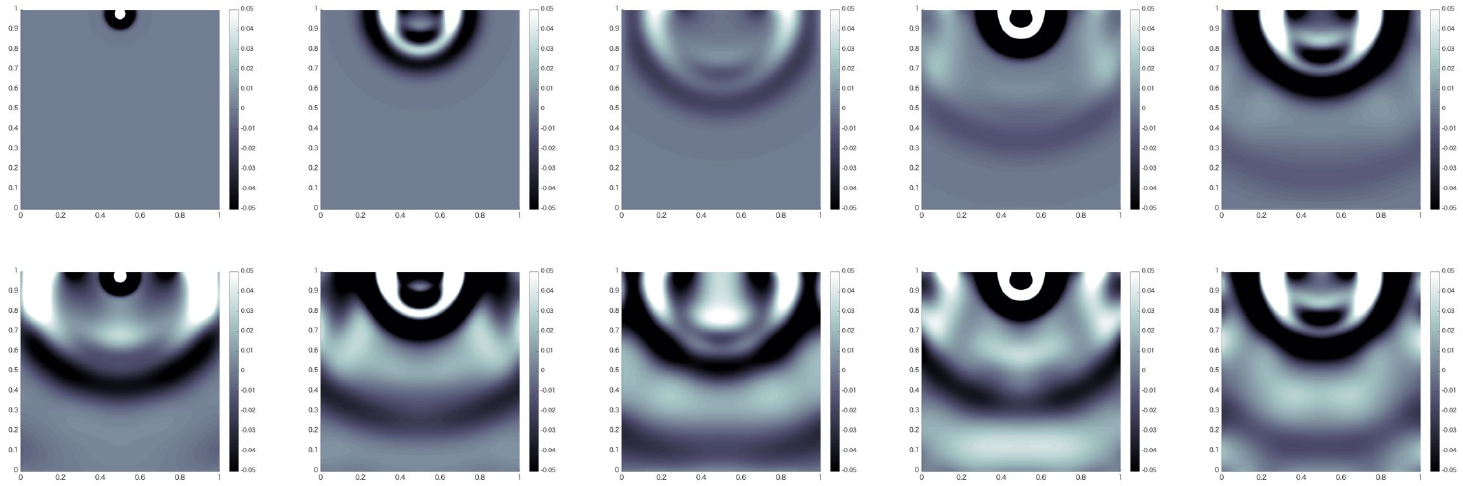}]
  \caption{Evolution of the dilation of $\hat{\u}$. Time is from top to bottom, left to right.}
  \label{v}
\end{figure}

\begin{figure}[H]
  \centering
  \includegraphics[page=3,scale=0.6]{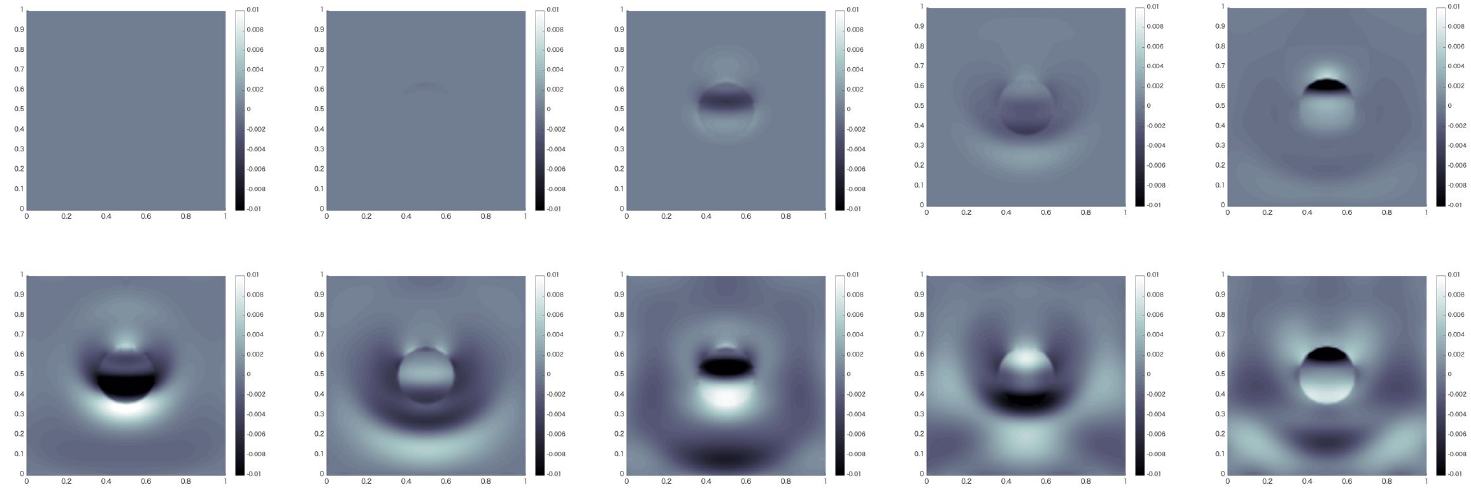}
  \caption{Difference in the dilations: $\u-\hat{\u}$.}
  \label{uv}
\end{figure}

For comparison, we perform the same simulation where $\Omega$ is composed of entirely $\Omega_2$. The result is shown in figure $\ref{v}$, where we again display the dilation $\hat{\u}$. The difference between the dilations are shown in figure $\ref{uv}$, where we clearly observe differences at the location of the inclusion.

\section{The inverse problem}\label{inv}
Having stated the model equation for elastic wave propagation through composite media and developed numerical methods for investigating its numerical simulation, we can now state the target inverse problem of this research. Moreover, for sake of agreement with the numerical investigations of the inverse problem, we will assume the case of isotropic materials. We remark that, in this setting, the model equation (\ref{eq:composite wave}) becomes the composite linear elastic wave equation (see e.g., \cite{Laude}).

As before, we assume that the domain $\Omega$ is composed of $\Omega_1$ and $\Omega_2$ and that their corresponding level set function has been constructed:
\begin{align}
  \theta(\x) =
  \begin{cases}
    -\text{dist}_{\Omega_1}(\x) & \x \in \Omega_1\\
    \text{dist}_{\Omega_2}(\x) & \x \in \Omega_2,\\
  \end{cases}
\end{align}
where $\text{dist}_{E}(\x)$ refers to the distance function (\ref{eq:dist}) to the set $E$.

For a given layout of $\Omega_1$ and $\Omega_2$, we define a target data $\delta$ (defined below) obtained from the solution of
the model equation:
\begin{align}
  \begin{cases}
    \rho(\x)\boldsymbol{u}_{tt} = \div \left( \sigma(\x) \right)+\f(\x,t)&\text{ in }(0,T)\times\Omega\\
    u_i = g_i &\text{ on }(0,T)\times\Gamma_{g_i} \text{ (Dirichlet)}\\
    \sigma_{ij}n_j=h_i &\text{ on }(0,T)\times\Gamma_{h_i} \text{ (Neumann)}\\
    u_i(0,\x)=u_{0i}(\x) &\x \in \Omega\\
    u_{i,t}(0,\x)=\dot{u}_{0i}(\x)& \x \in \Omega\\
    [u_i]=0 & \text{on}\;\partial\Omega_1\cap\partial\Omega_2\hspace{15pt} t\in(0,T)\\
    [\sigma_{ij}n_{j}]=0 & \text{on}\;\partial\Omega_1\cap\partial\Omega_2 \hspace{15pt} t\in(0,T)
  \end{cases}
  \label{eq:star}
\end{align}
for given initial and boundary conditions, and where $\sigma(\x)$ denotes the spatially dependent stress tensor:
\begin{align}
  \sigma(\x) = c(\x):\epsilon(\u),
\end{align}
where $\epsilon(\u)$ denotes the strain tensor:
\begin{align}
  \epsilon(\u) = \frac{1}{2}\left( \nabla \u + (\nabla \u)^{T}\right).
\end{align}
Here we remark that the density and stress tensor functions can be controlled through the level set function $\theta(\x)$. \noindent
The interfacial boundary condition $[u_i]=0$ specifies continuity in the displacement field, whereas $[\sigma_{ij}n_{j}]=0$ designates continuity of the stresses normal to the interface (see e.g., \cite{Laude}). Here $n_j$ denotes the $j^{th}$ component of the normal to the interface, $\n$, and we assume that the outer force $\f(\x,t)$ is a given time dependent vector field, such as equation (\ref{mylaser}).

\begin{figure}[H]
  \centering
  \includegraphics[scale=0.35]{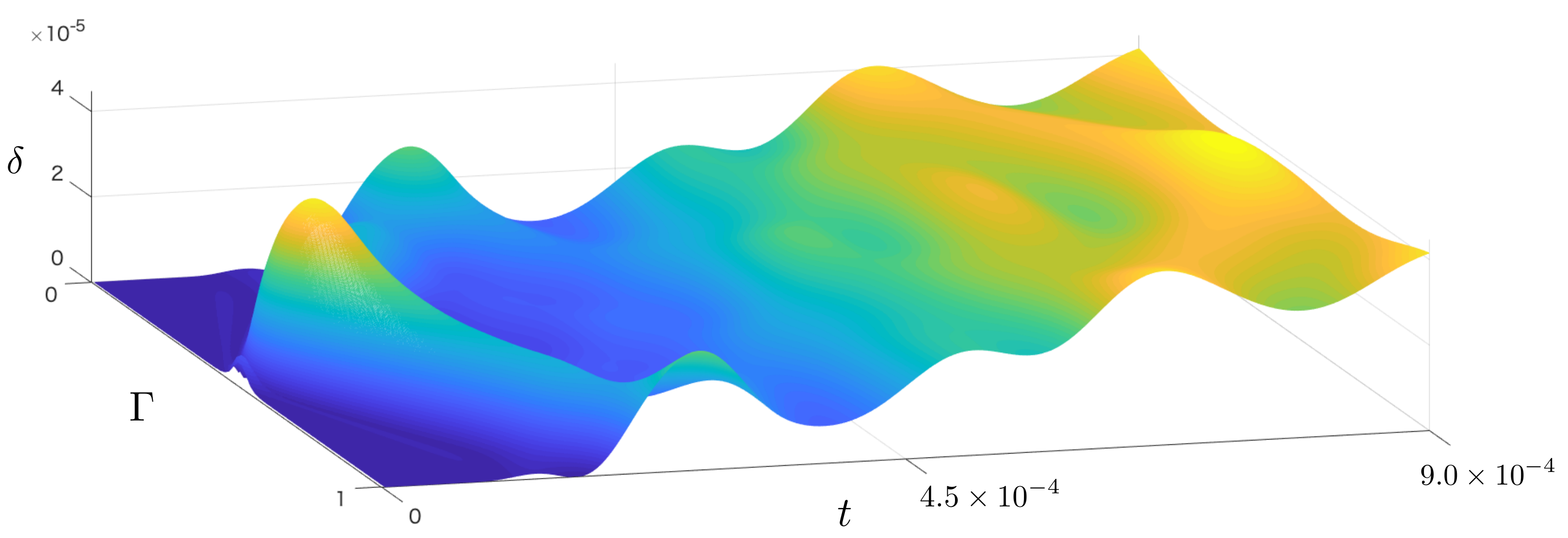}
  \caption{The boundary data $\delta$ on $\Gamma$.}
  \label{delta}
\end{figure}

Given a solution to the model equation ($\ref{eq:star}$), the boundary data $\delta$ is defined as a single component of the displacement
field along an external boundary $\Gamma$ of $\Omega$:
\begin{align}
\delta(\y,t) = u_{i}(\y,t) \hspace{15pt} (\y,t)\in\Gamma\times (0,T).
\end{align}
For example, when the closure of $\Omega$ is the closed unit cube $[0,1]\times[0,1]$, we can define $\Gamma:= \{\x\in \Omega\hspace{2pt}|\hspace{2pt} y=1\}$ and $\delta(x,t)=u_2(x,1,t).$ Figure \ref{delta} shows the evolution of $\delta$ on the time interval $[0,9.0\times 10^{-4}$, where $\Omega_1$, $\Omega_2$, as well as their boundary conditions and material parameters are prescribed as in the simulations of the  previous section (see equation (\ref{eq:location})). Here we have prescribed $\f(\x,t)$ using equation (\ref{mylaser}). We now state the target inverse problem of this study.

Given a boundary data $\delta: \Gamma\times(0,T)\rightarrow {\bf{R}}$ and a penalty parameter $\tilde{\epsilon}>0$, find $\theta(\x)$ minimizing the cost functional
\begin{equation}
  \mathcal{E}(\u,\delta)=\frac{1}{\tilde{\varepsilon}}||u_2 - \delta||^2_{L^2(\Gamma\times(0,T))},
  \label{eq:cost}
\end{equation}
where $\u(\x,t)$ solves the model equation:
\begin{align}
  \begin{cases}
    \rho\boldsymbol{u}_{tt} = \div \left( \sigma \right)+\f(\x,t)&\text{ in }(0,T)\times\Omega\\
    u_i = g_i &\text{ on }(0,T)\times\Gamma_{g_i} \text{ (Dirichlet)}\\
    \sigma_{ij}n_j=h_i &\text{ on }(0,T)\times\Gamma_{h_i} \text{ (Neumann)}\\
    u_i(0,\x)=u_{0i}(\x) &\x \in \Omega\\
    u_{i,t}(0,\x)=\dot{u}_{0i}(\x)& \x \in \Omega\\
    [u_i]=0 & \text{on}\;\partial\Omega_1\cap\partial\Omega_2\hspace{15pt} t\in(0,T)\\
    [\sigma_{ij}n_{j}]=0 &\text{on}\;\partial\Omega_1\cap\partial\Omega_2\hspace{15pt} t\in(0,T).
  \end{cases}
  \label{eq:inv}
\end{align}
Here the level set function $\theta({\x})$ controls the composition of $\Omega$, including the density function and elastic coefficients:
\begin{eqnarray}
  \Omega_{1} &=&\{\x\in \Omega | \theta(x) < 0\}\nonumber\\
  \Omega_{2} &=& \{\x\in \Omega | \theta(x) \geq 0\}\nonumber\\
  \rho(\theta(\x)) &=& \phi(\theta(\x))\rho_{1} + (1-\phi(\theta(\x)))\rho_{2}\nonumber\\
  c(\theta(\x))  &=&  \phi(\theta(\x))c_{1} + (1-\phi(\theta(\x)))c_{2}\nonumber\\
  \phi(\theta(\x)) &=& \dfrac{1}{2}\left(\tanh\left(\dfrac{\theta(\x)}{\varepsilon}\right)+1\right).\nonumber
\end{eqnarray}

\section{The method of Lagrange multipliers for the inverse problem}\label{lag}
In this section we derive the gradient of the cost functional $\mathcal{E}$ (see equation (\ref{eq:cost})) for investigating solutions of the target inverse problem (\ref{eq:inv}). In particular, upon determining this gradient, we will investigate the solution of the inverse problem as the gradient flow of the cost functional. In order to do this, we will use a level set function $\theta(\x)$ to describe the shape and location of the composite materials. This function will serve as the control parameter for minimizing $\mathcal{E}$. We will now employ the method of Lagrange multipliers, suitably adapted to the functional setting (see, e.g., \cite{Azegami}, \cite{Laporte}).
\subsection{Gradient of the cost functional $\mathcal{E}$}

In order to use the method of Lagrange multipliers, we define the following Lagrangian:
\begin{equation}
  \mathcal{L}(\theta, \u,\v) = \mathcal{E}(\theta,\u) + \mathcal{L}_{S}(\theta, \u,\v),
\end{equation}
where $\u$ and $\v$ are arbitrary vector valued functions and $\theta(\x)$ denotes the level set function (\ref{eq:levelset}). Here, 
$\mathcal{E}$ denotes the cost functional (\ref{eq:cost}) for a given penalty parameter $\tilde{\epsilon}>0:$
\begin{equation}
  \mathcal{E}(\u,\delta)=\frac{1}{\tilde{\varepsilon}}||u_2 - \delta||^2_{L^2(\Gamma\times(0,T))}
\end{equation}
and $\mathcal{L}_{S}$ is expressed by the following equation:
\begin{equation}
  \mathcal{L}_{S}(\theta, \u,\v) = \int_{0}^{T}\int_{\Omega}\left( -\rho(\theta)\u_t\cdot \v_t + (c(\theta)\otimes\nabla \u):\nabla\v + \f\cdot \v\right){dxdt}.
\end{equation}
In the above, $T$ is a positive value corresponding to the final time in ($\ref{eq:composite wave}$) and the vector valued function $\f(\x,t)$ represents an outer force (for example, given by (\ref{mylaser})). For $A,B\in\R^{n\times n}$, the symbol $:$ in  $(c(\theta)\otimes\nabla \u):\nabla\v$ refers to the quantity:
\begin{equation}
  A:B = \sum_{j=1}^n\sum_{i = 1}^{n} A_{ji}B_{ji}.
\end{equation}

The Frechet derivative of the Lagrangian $\mathcal{L}$ is expressed:
\begin{equation}
  \mathcal{L}'(\theta, \u,\v) =
  \mathcal{L}_{\theta}(\theta, \u,\v)[\vartheta] +
  \mathcal{L}_{\u}(\theta, \u,\v)[\tilde{\u}]+
  \mathcal{L}_{\v}(\theta, \u,\v)[\tilde{\v}],
  \label{eq:total}
\end{equation}
where $\vartheta, \tilde{\u},$ and $\tilde{\v}$ are arbitrary variations of $\theta, \u$, and $\v$, respectively.

The partial Frechet derivative of $\mathcal{L}$ with respect to $\v$ is expressed:
\begin{align}
  \begin{aligned}
   &\mathcal{L}_{\v}(\theta, \u,\v)[\tilde{\v}]\\
    &=\lim_{\varepsilon=0}\dfrac{\partial}{\partial\epsilon}
    \int_{0}^{T}\int_{\Omega}\left( -\rho(\theta)\u_t\cdot (\v
    +\epsilon\tilde{\v})_t
    +(c(\theta)\otimes\nabla \u):\nabla(\v+\epsilon\tilde{\v})
    +\f\cdot(\v+\epsilon\tilde{\v})\right){dxdt}.\nonumber\\
  \end{aligned}
\end{align}
A formal computation then yields:
\begin{align}
  \begin{aligned}
    \mathcal{L}_{\v}(\theta, \u,\v)[\tilde{\v}]
    =\int_{0}^{T}\int_{\Omega}(-\rho(\theta)\u_{t}\cdot\tilde{\v}_t
    +(c(\theta)\otimes\nabla \u):\nabla \tilde{\v} +\f\cdot\tilde{\v}){dxdt}
    =\mathcal{L}_S(\theta, \u,\tilde{\v}).\nonumber
  \end{aligned}
\end{align}
When $\u$ is a solution of equation ($\ref{eq:composite wave}$), we remark that $\mathcal{L}_S(\theta, \u,\v)$ takes the value 0 for
any arbitrary $\v$.

Similarly, the partial Frechet derivative of $\mathcal{L}$ with respect to $\u$ is computed:
\begin{align}
  \begin{aligned}
    \mathcal{L}_{\u}(\theta, \u,\v)[\tilde{\u}]
    =&\lim_{\varepsilon=0}\dfrac{\partial}{\partial \epsilon}\int_{0}^{T}\int_{\Gamma} \frac{1}{\tilde{\varepsilon}}((\u+\epsilon\tilde{\u})_2-\delta)^2 dSdt\\
    &+\int_{0}^{T}\int_{\Omega}-\rho(\theta)(\u+\epsilon\tilde{\u})_t\cdot\v_t
    +(c(\theta)\otimes\nabla (\u+\epsilon\tilde{\u})) : \nabla \v {dxdt}\\
    =&\int_{0}^{T}\int_{\Gamma}\frac{2}{\tilde{\varepsilon}}(u_2-\delta)\tilde{u}_2dSdt\\
    \;&+\int_{0}^{T}\int_{\Omega}-\rho(\theta)\tilde{\u}_t\cdot\v_t
    +\left(c(\theta)\otimes\nabla \tilde{\u}\right) : \nabla \v {dxdt}.
  \end{aligned}
\end{align}
Here we remark that $\mathcal{L}_{\u}(\theta,\u,\v)[\tilde{\u}]$ takes the value 0 whenever $\v$ is a solution of the following adjoint problem:
\begin{align}
  \begin{cases}
    \rho(\theta)\boldsymbol{v}_{tt} = \div \left( c(\theta)\otimes\nabla\v \right)&\text{ in }(0,T)\times\Omega\\
    v_2 = \dfrac{2}{\tilde{\epsilon}}(u_2 - \delta) &\text{ on }(0,T)\times\Gamma\\
    v_i = g_i &\text{ on }(0,T)\times\Gamma_{g_i} \text{ (Dirichlet)}\\
    (c(\x)\otimes\nabla\v)_{ij} n_j=0 &\text{ on }(0,T)\times\Gamma_{h_i} \text{ (Neumann)}\\
    v_i(0,\x)=u_{0i}(\x) &\x \in \Omega\\
    v_{i,t}(0,\x)=\dot{u}_{0i}(\x)& \x \in \Omega\\
    [v_i]=0 & \text{on}\;\partial\Omega_1\cap\partial\Omega_2\hspace{15pt} t\in(0,T)\\
    [(c(\theta)\otimes\nabla\v)_{ij}n_{j}]=0 & \text{on}\;\partial\Omega_1\cap\partial\Omega_2\hspace{15pt} t\in(0,T).
  \end{cases}
  \label{eq:myadjoint}
\end{align}
In this way, having set $\u$, and $\v$, and since $\vartheta$ is an arbitrary variation of $\theta$, we can obtain the gradient of the cost functional as the remaining term in equation (\ref{eq:total}):
\begin{align}
  \begin{aligned}
    \mathcal{L}'(\theta, \u,\v)
    &=\mathcal{L}_{\theta}(\theta, \u,\v)[\vartheta]\\
    &=\int_{0}^{T}\int_{\Omega} (-\rho'(\theta)\u_t\cdot\v_t
    +\left(c'(\theta)\otimes\nabla \u\right):\nabla \v)\vartheta {dxdt}\\
    &=\langle g,\vartheta\rangle_{L^2(\Omega)}.
  \end{aligned}
\end{align}

Here $\rho'(\theta),c'(\theta)$ are regularizations of delta functions supported on interfaces between the regions $\Omega_1$ and $\Omega_2$:
\begin{align}
  \rho'(\theta(\x)) &= (\rho_2-\rho_1)\phi'(\theta(\x))\\
  c'(\theta(\x)) &= (c_2-c_1)\phi'(\theta(\x))\\
   \phi'(\theta(\x)) &= \frac{1}{2\varepsilon}\text{sech}^2\left(\dfrac{\theta(\x)}{\varepsilon}\right).
\end{align}
In the above, $\varepsilon$ is a positive parameter which prescribes the cost functional's gradient to act as a regularized delta function near the interface. We remark that the particular value of $\varepsilon$ is important when performing numerical calculations.

Having determined the gradient of the cost functional, we now turn to designing an approximation method for solving the
inverse problem.

%
%

%
%
%
%
\section{An approximation method for the inverse problem}\label{level}
In this section, using the cost functional's gradient information obtain by the method of Lagrange multipliers from the previous section, we will develop an approximation scheme for investigating the numerical behavior of solutions to our target inverse problem. A brief explanation of our method is as follows. We first construct the boundary data $\delta$ for a given target composition. Second, we set the composition of the region to be different from the target and solve the corresponding model equation. This gives us the displacement field arising in the adjoint problem (\ref{eq:myadjoint}). After solving the adjoint problem we are then able to construct the gradient of the cost functional. This gradient information is then used within the level set method to evolve the interface and decrease the value of the cost functional. The details are explained in the following.

\subsection{Algorithm for approximating solutions of the inverse problem}
The numerical solution of the inverse problem is expressed as a gradient flow of the cost functional (\ref{eq:cost}), which is performed as follows.

\noindent \underline{Step 1: Construct the boundary data $\delta$.}
\begin{enumerate}
  \item Set the target composition of $\Omega$ and construct its level set function:
  \begin{align}
    \theta_{target}(\x) =
    \begin{cases}
      -\text{dist}_{\tilde{\Omega}_1}(\x) & \x \in \tilde{\Omega}_1\\
      \text{dist}_{\tilde{\Omega}_2}(\x) & \x \in \tilde{\Omega}_2,
    \end{cases}
  \end{align}
  where $\Omega = \tilde{\Omega}_1\cup\tilde{\Omega}_2$ denotes the {\emph{target}} domain composition.
  \item Solve the model equation (\ref{eq:optim1}) corresponding to $\theta_{target}(\x)$:
    \begin{align}
    \begin{cases}
      \rho(\theta_{target}(\x))\boldsymbol{u}_{tt} = \div \left( \sigma(\theta_{target}(\x)) \right)+\f(\x,t)&\text{ in }(0,T)\times\Omega\\
      u_i = g_i &\text{ on }(0,T)\times\Gamma_{g_i} \text{ (Dirichlet)}\\
      \sigma_{ij}(\theta_{target}(\x))n_j=h_i &\text{ on }(0,T)\times\Gamma_{h_i} \text{ (Neumann)}\\
      u_i(0,\x)=u_{0i}(\x) &\x \in \Omega\\
      u_{i,t}(0,\x)=\dot{u}_{0i}(\x)& \x \in \Omega\\
      [u_i]=0 & \text{on}\;\partial\tilde{\Omega}_1\cap\partial\tilde{\Omega}_2\hspace{15pt} t\in(0,T)\\
      [\sigma_{ij}(\theta_{target}(\x))n_{j}]=0 & \text{on}\;\partial\tilde{\Omega}_1\cap\partial\tilde{\Omega}_2\hspace{15pt} t\in(0,T),
    \end{cases}
    \label{eq:optim1}
  \end{align}
  for a given outer force $\f(\x,t)$, and initial and boundary conditions.
  \item Extract the boundary data from the solution:
  \begin{align}
    \delta(\y,t) = \{u_2(\y,t)\hspace{2pt}|\hspace{2pt}\y \in \Gamma ,t>0\},
  \end{align}
  where $\Gamma$ denotes an exterior boundary of $\Omega$.
\end{enumerate}
\vspace{15pt}
\noindent \underline{Step 2: Set the initial condition for the level set method.}
\begin{align}
  \theta_{init}(\x) =
  \begin{cases}
    -\text{dist}_{\Omega_1}(\x) & \x \in \Omega_1\\
    \text{dist}_{\Omega_2}(\x) & \x \in \Omega_2\\
  \end{cases}
\end{align}
Here, $\Omega_1$ and $\Omega_2$ are generally assumed to be different from the target compositions $\tilde{\Omega}_1$ and 
$\tilde{\Omega}_2$.
\vspace{15pt}

\noindent \underline{Step 3: Construct the cost functional gradient of $\mathcal{E}$ with respect to $\theta(\x)$.}
\begin{enumerate}
  \item Solve the model equation (\ref{eq:optim2}) corresponding to $\theta_{init}(\x)$:
  \begin{align}
    \begin{cases}
      \rho(\theta_{init}(\x))\boldsymbol{u}_{tt} = \div \left( \sigma(\theta_{init}(\x)) \right)+\f(\x,t)&\text{ in }(0,T)\times\Omega\\
      u_i = g_i &\text{ on }(0,T)\times\Gamma_{g_i} \text{ (Dirichlet)}\\
      \sigma_{ij}(\theta_{init}(\x))n_j=h_i &\text{ on }(0,T)\times\Gamma_{h_i} \text{ (Neumann)}\\
      u_i(0,\x)=u_{0i}(\x) &\x \in \Omega\\
      u_{i,t}(0,\x)=\dot{u}_{0i}(\x)& \x \in \Omega\\
      [u_i]=0 & \text{on}\;\partial\Omega_1\cap\partial\Omega_2\hspace{15pt} t\in(0,T)\\
      [\sigma_{ij}(\theta_{init}(\x))n_{j}]=0 &\text{on}\;\partial\Omega_1\cap\partial\Omega_2\hspace{15pt} t\in(0,T)
    \end{cases}
    \label{eq:optim2}
  \end{align}
  \item Solve the adjoint problem $(\ref{eq:optim3})$, using $u_2(\x,t)$ and $\delta$:
  \begin{align}
    \begin{cases}
      \rho(\theta_{init}(\x))\boldsymbol{v}_{tt} = \div \left( \sigma(\theta_{init}(\x)) \right)&\text{ in }(0,T)\times\Omega\\
      v_2 = \dfrac{2}{\tilde{\epsilon}}(u_2 - \delta) &\text{ on }(0,T)\times\Gamma\\
      v_i = g_i &\text{ on }(0,T)\times\Gamma_{g_i} \text{ (Dirichlet)}\\
      \sigma_{ij}(\theta_{init}(\x))n_j=h_i &\text{ on }(0,T)\times\Gamma_{h_i} \text{ (Neumann)}\\
      v_i(0,\x)=u_{0i}(\x) &\x \in \Omega\\
      v_{i,t}(0,\x)=\dot{u}_{0i}(\x)& \x \in \Omega\\
      [v_i]=0 & \text{on}\;\partial\Omega_1\cap\partial\Omega_2\hspace{15pt} t\in(0,T)\\
      [\sigma_{ij}(\theta_{init}(\x))n_{j}]=0 &\text{on}\;\partial\Omega_1\cap\partial\Omega_2\hspace{15pt} t\in(0,T)
    \end{cases}
    \label{eq:optim3}
  \end{align}
  \item Construct the cost functional gradient of $\mathcal{E}$ as the function:
  \begin{align}
    \begin{aligned}
      g(\x)
      &=\int_{0}^{T}(-\rho'(\theta_{init})\u_t\cdot\v_t
      +\left(c'(\theta_{init})\otimes\nabla \u\right) : \nabla \v){dt}
      \label{eq:optim4}
    \end{aligned}
  \end{align}
  where $\rho(\theta_{init})$ and $c(\theta_{init})$ are given by equation ($\ref{eq:parameters}$).
\end{enumerate}
\vspace{15pt}
\noindent \underline{Step 4: Evolve the level set function using the gradient $g(\x)$ within the level set equation.}
\begin{enumerate}
  \item Solve the level set equation for a time $\tau^*>0 $:
  \begin{align}
    \begin{cases}
      &\theta_{\tau}={|\nabla\theta|}\displaystyle\int_{0}^{T}\rho'(\theta)\u_t \cdot \v_t + (c'(\theta)\otimes \nabla\u):\nabla \v{dt} \hspace{15pt} \tau\in(0,\tau^*),\x\in\Omega)\\
      &\theta(\x,0)=\theta_{init}(\x)  \hspace{147pt} \x \in \Omega
    \end{cases}
    \label{eq:level}
  \end{align}
  \item Evolve the level set function:
  \begin{align}
    \theta_{init}(\x) =
    \begin{cases}
      -\text{dist}_{\theta(\tau^*)>0}(\x) & \x\in\theta(\tau^*)>0\\
      \text{dist}_{\theta(\tau^*)\le 0}(\x) & \x\in\theta(\tau^*)\le 0
    \end{cases}
  \end{align}
\end{enumerate}
\vspace{15pt}
\noindent \underline{Step 5: Repeat steps 3 and 4 until $\mathcal{E}$ converges.}

\subsection{Properties of the cost functional gradient}\label{propgrad}
As we have explained, the gradient of the cost functional depends on various parameters through solutions to partial differential equations. For this reason, we have performed numerical experiments to confirm its basic properties. In particular, in this section, we will show that the cost functional gradient indeed holds information on the location of inclusions. In particular, using the domain triangulation and parameters set forth in section \ref{numbehave},  we set set the initial condition for the level set method as in equation (\ref{eq:location}) and investigate the behavior of $g(\x)$ for seven different target compositions. 

We employ an immersed boundary technique in the adjoint equation for $\v$. In particular, by setting Neumann zero boundary conditions on the boundary of $\Omega$, our approach models the Dirichlet boundary condition on $\Gamma$ by using a regularized delta function as an outer force:
\begin{align}
\f_{\Gamma}(\x,t)=\frac{2}{\tilde{\epsilon}}\left(u_2(\x,t) - \delta(x,t) \right)\frac{\text{sech}^{2}(y-l_y)}{2\epsilon},
\end{align}
where we set $\tilde{\epsilon}=1/1000$, and $\epsilon=1/25$.

\begin{figure}[H]
  \centering
  \includegraphics[scale=0.5]{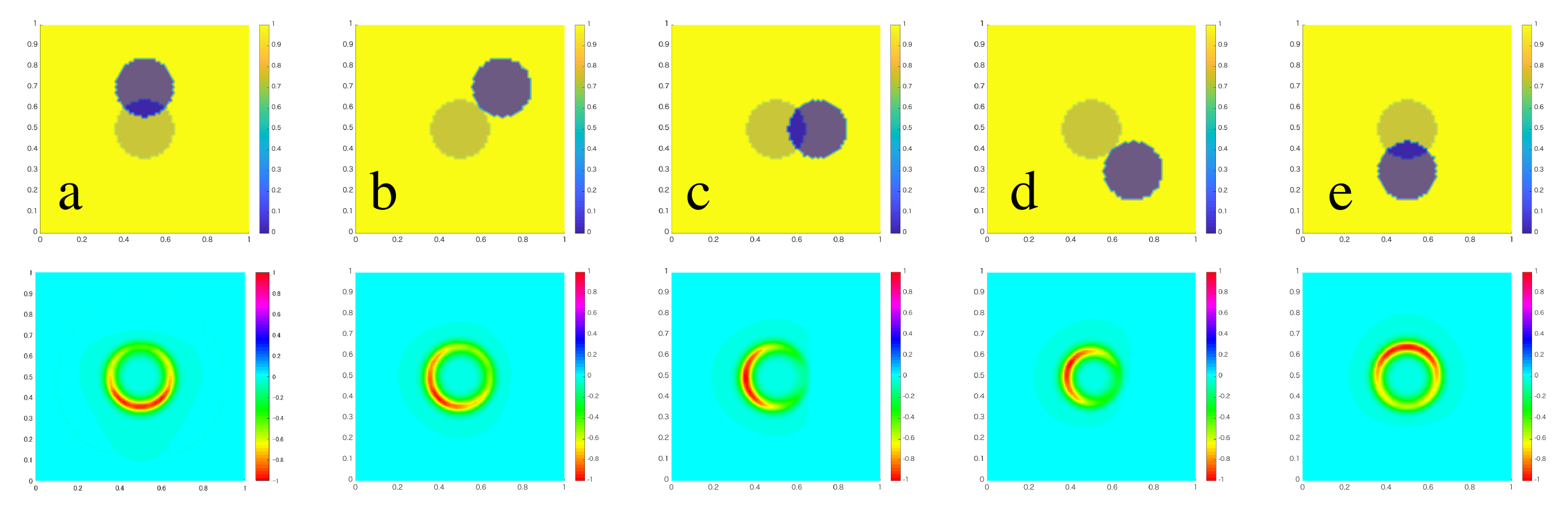}
  \caption{Top row: $\theta_{target}$ and $\theta_{init}$. Bottom row: the corresponding gradients.}
  \label{fig:g_confirm}
\end{figure}

Figure \ref{fig:g_confirm} shows the cost functional gradient $g$ responding to various target compositions. The first row of the figure shows the initial condition (the light blue region) overlaid with each target (shown in dark blue). The second row shows the corresponding gradients.

\begin{figure}[H]
  \centering
  \includegraphics[scale=0.5,page=2]{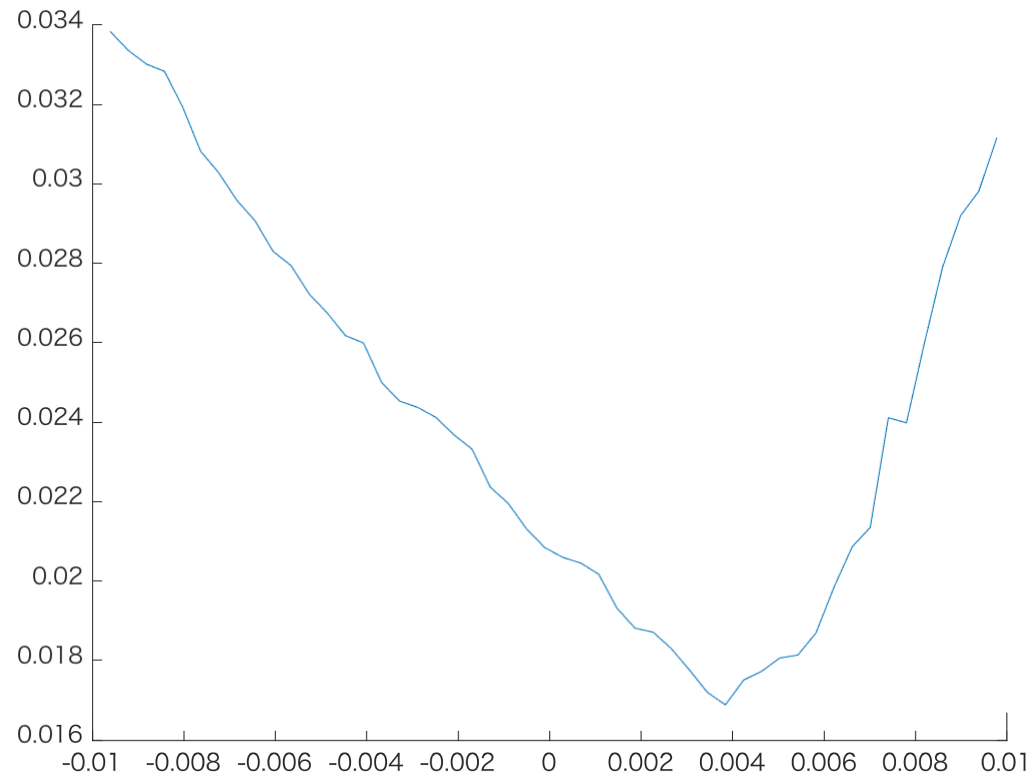}
  \caption{Cost functional values: $\mathcal{E}(\theta_{init}(\x) - \tau g(\x)),\;\tau \in [-0.01,0.01]$.}
  \label{fig:tau5}
\end{figure}

Six cases were considered. In the first (figure \ref{fig:g_confirm} a), the $\theta_{target}$ is placed above
$\theta_{init}$ and near the upper boundary. Inspection of the gradient $g$ confirms that the level set equation evolves the interface with a normal velocity in the direction of the target. As will be shown in the next section, evolving the level set function by the level set equation decreases both the distance between the target and the value of the cost functional. Figure \ref{fig:tau5} displays the value of the cost functional $\mathcal{E}(\theta_{init}-\tau g)$ for $\tau$ within the interval $[-0.01, 0.01]$ where we are able to verify that the gradient descent of $\mathcal{E}$ is in the direction $-g(\x)$.


\section{Numerical solution of the inverse problem}\label{numresult}
This section examines the numerical solution of the inverse problem by observing the gradient flow of the cost functional. We begin by employing our approximation method on each of the five configurations shown in figure \ref{fig:g_confirm} a-e.

\begin{figure}[H]
  \centering
  \includegraphics[scale=0.5,page=2]{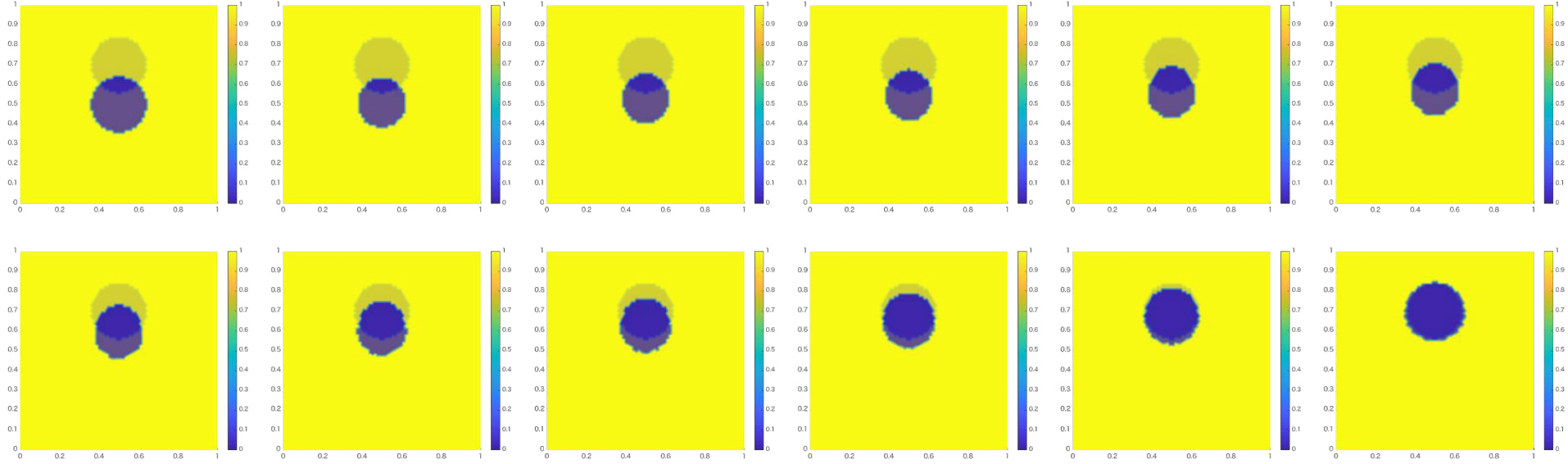}
  \caption{Numerical solution of the inverse problem. Time is from top to bottom, left to right. The dark region is the approximate solution and $\theta_{target}$ is shown as the light disk.)}
  \label{fig:optim5}
\end{figure}

\begin{figure}[H]
  \centering
  \includegraphics[scale=0.4,page=1]{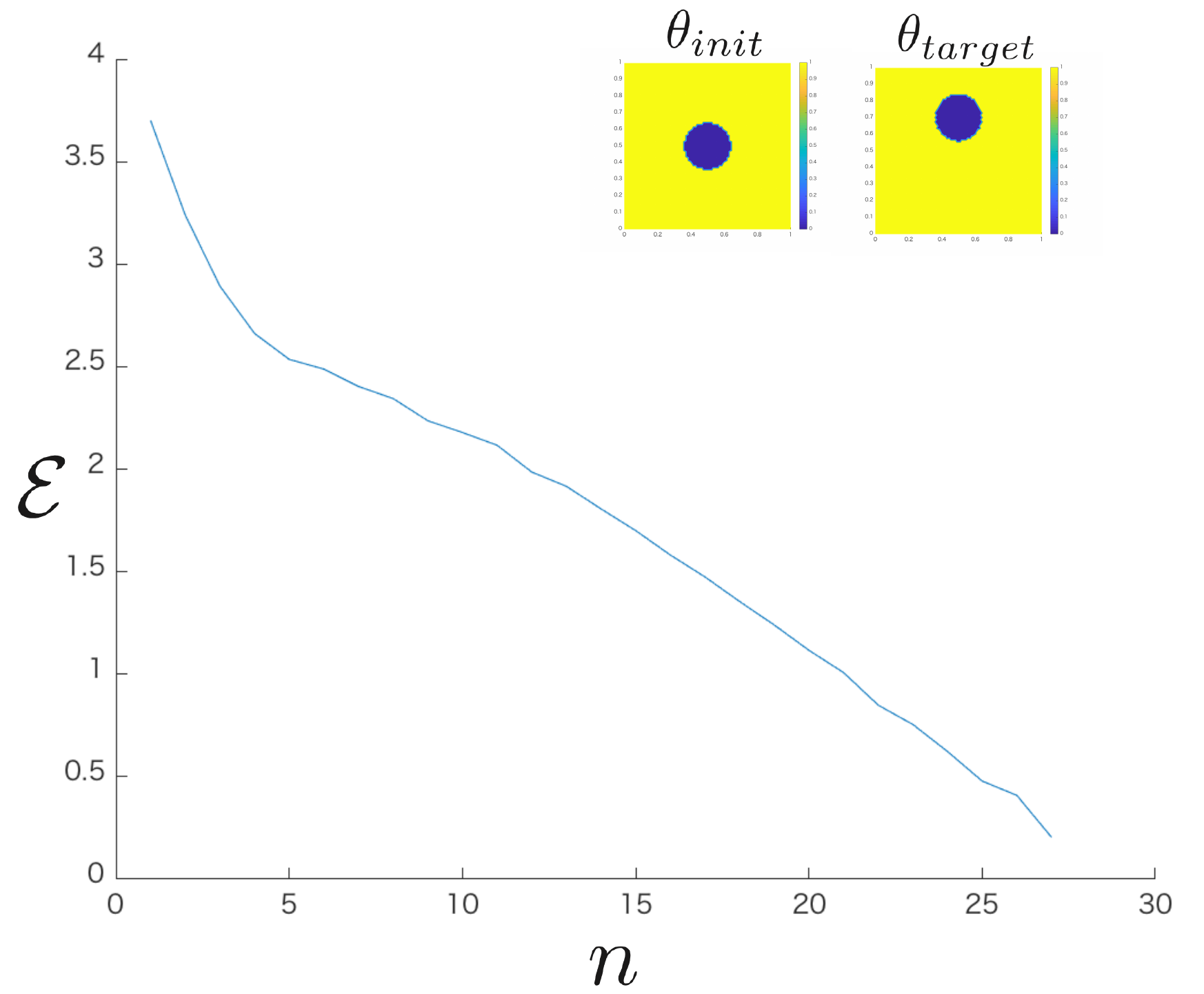}
  \caption{Evolution of the cost functional $\mathcal{E}$ under the gradient flow. The initial and target conditions are also displayed.}
  \label{fig:graph2}
\end{figure}

The level set function $\theta_{init}(\x)$, the outer force $\f(\x,t)$, and the remaining parameters are set as in the previous sections. We set $\theta_{target}(\x)$ (a disk with the same radius as $\theta_{init}(\x)$, shown as a light shadow in the figures) directly above the initial condition and near the boundary of the computational domain (see figure \ref{fig:g_confirm} a). Repeatedly solving the level set equation (\ref{eq:level}) yields the evolution shown in figure \ref{fig:optim5}. We have displayed twelve instants of the evolution, where we observe the level set equation's solutions to converge to the prescribed target.  Figure \ref{fig:graph2} shows a plot of the cost functional's value as a function of the $n^{th}$ iteration of the approximation method. We confirm the value of the cost functional to decrease along the flow.

\begin{figure}[H]
  \centering
  \includegraphics[scale=0.5,page=1]{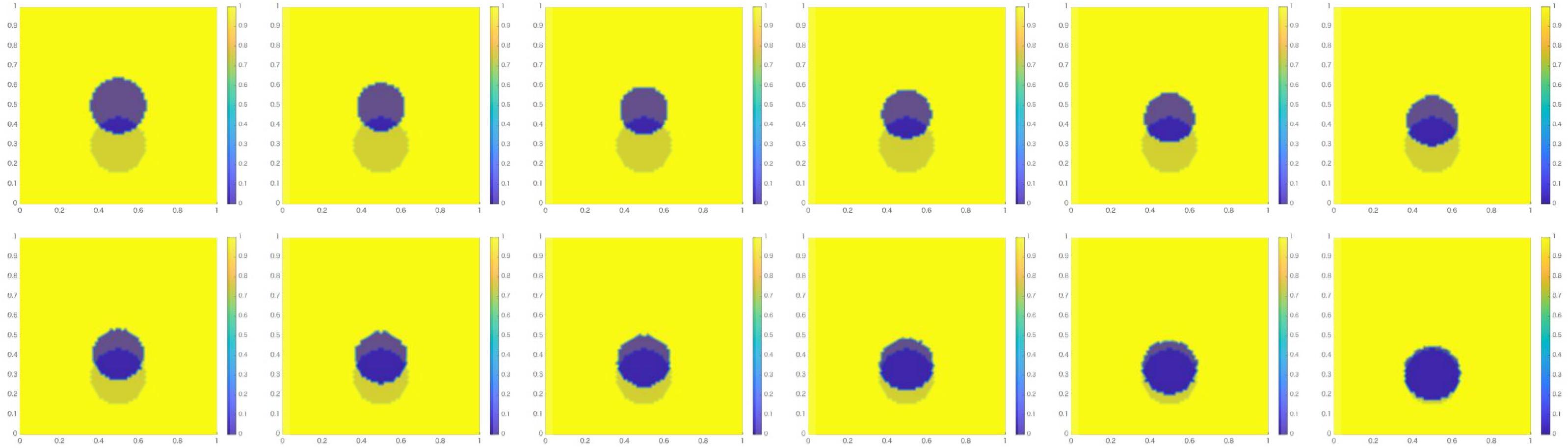}
  \caption{Numerical solution of the inverse problem. Time is from top to bottom, left to right. The dark region is the approximate solution and $\theta_{target}$ is shown as the light disk.)}
  \label{fig:optim4}
\end{figure}

\begin{figure}[H]
  \centering
  \includegraphics[scale=0.4,page=2]{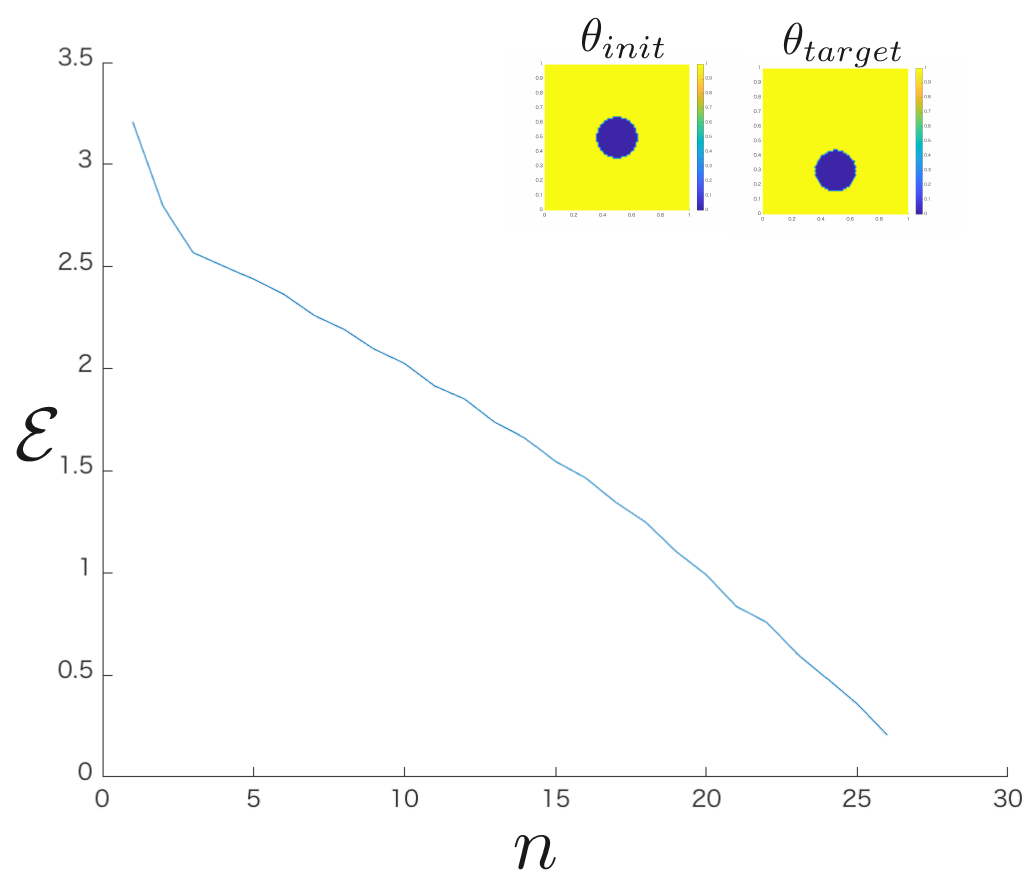}
  \caption{Evolution of the cost functional $\mathcal{E}$ under the gradient flow. The initial and target conditions are also displayed.}
  \label{fig:graph1}
\end{figure}

Our next computation sets $\theta_{target}(\x)$ directly below the initial condition, near the bottom boundary (see figure {\ref{fig:g_confirm}} e). We apply our approximation method and observe the motion of the level set function as it converges to its target. Figure \ref{fig:optim4} shows 12 instants in time of the evolution, and figure \ref{fig:graph1} shows the value of the cost functional at each iteration step. We again confirm that the gradient flow decreases the value of the cost functional.

\noindent

\begin{figure}[H]
  \centering
  \includegraphics[scale=0.5,page=3]{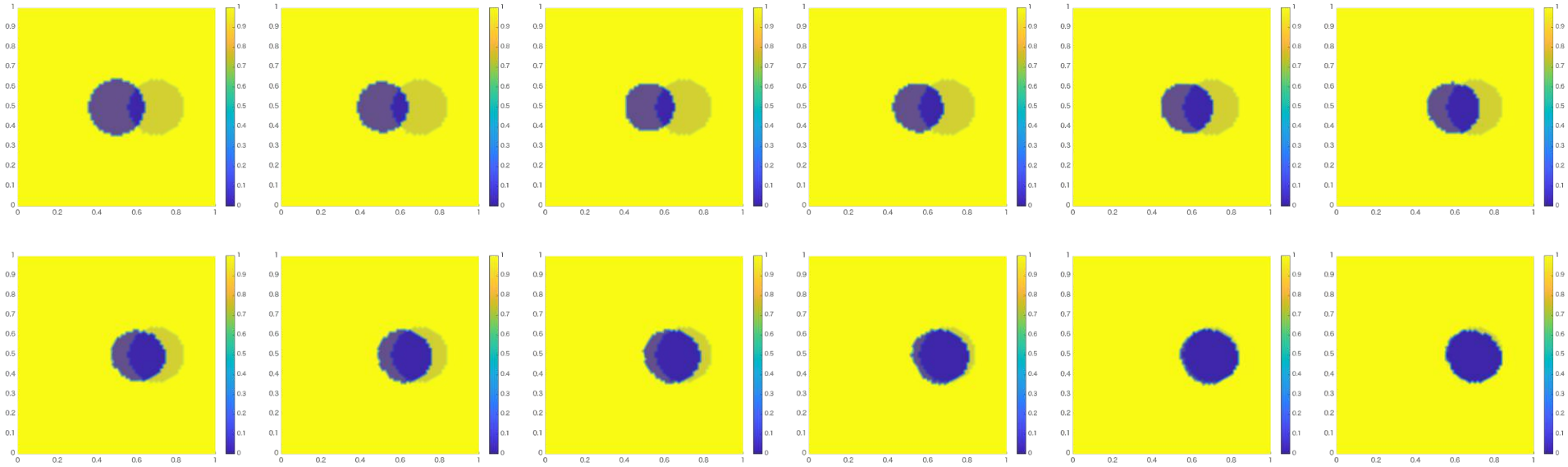}
  \caption{Numerical solution of the inverse problem. Time is from top to bottom, left to right. The dark region is the approximate solution and $\theta_{target}$ is shown as the light disk.)}
  \label{fig:optim6}
\end{figure}

\begin{figure}[H]
  \centering
  \includegraphics[scale=0.4,page=3]{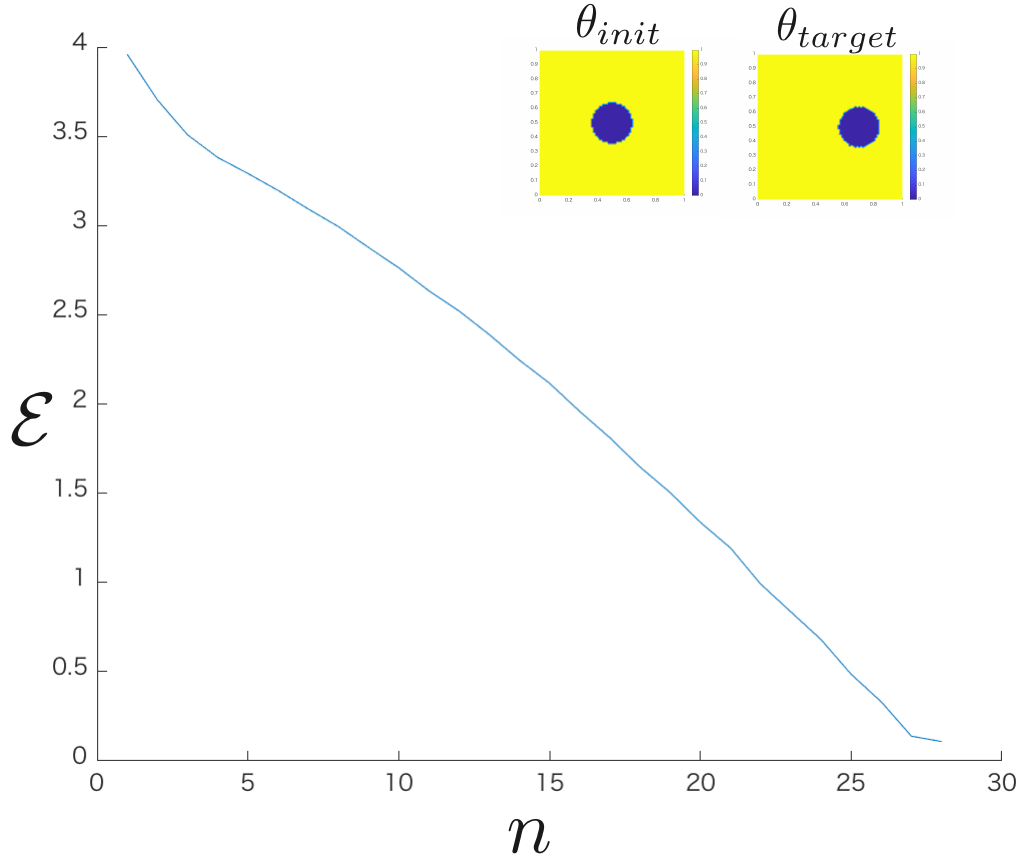}
  \caption{Evolution of the cost functional $\mathcal{E}$ under the gradient flow. Initial and target conditions are also displayed.}
  \label{fig:graph3}
\end{figure}

\noindent
A similar result can be obtained when the target is located to the right of the initial condition (see figure \ref{fig:g_confirm} c). Our method yields an approximation of the gradient flow, shown in figure \ref{fig:optim6}. We again observe the values of the cost functional to decrease with each iteration of the optimization (figure \ref{fig:graph3}).


\begin{figure}[H]
  \centering
  \includegraphics[scale=0.5,page=4]{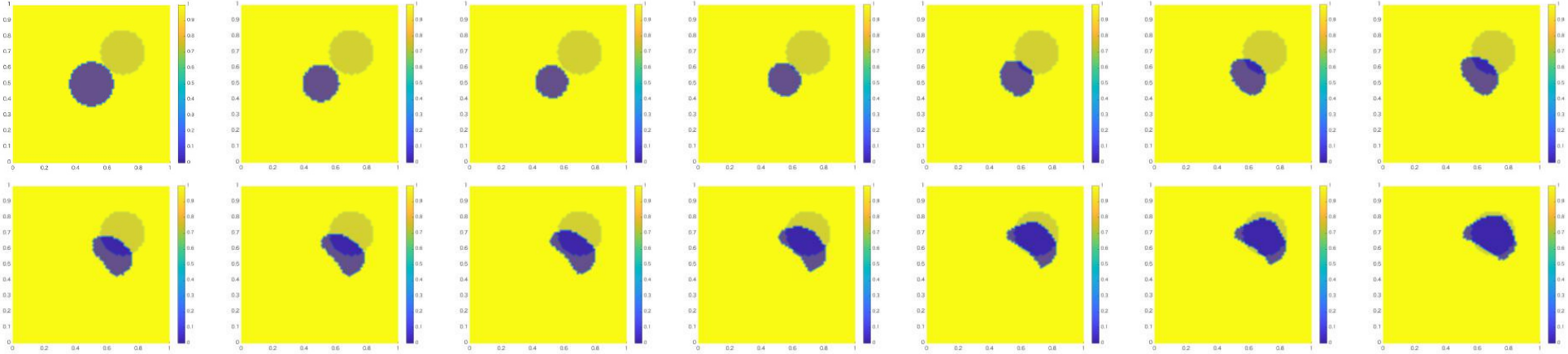}
  \caption{Numerical solution of the inverse problem. Time is from top to bottom, left to right. The dark region is the approximate solution and $\theta_{target}$ is shown as the light disk.}
  \label{fig:optim7}
\end{figure}

\begin{figure}[H]
  \centering
  \includegraphics[scale=0.4,page=4]{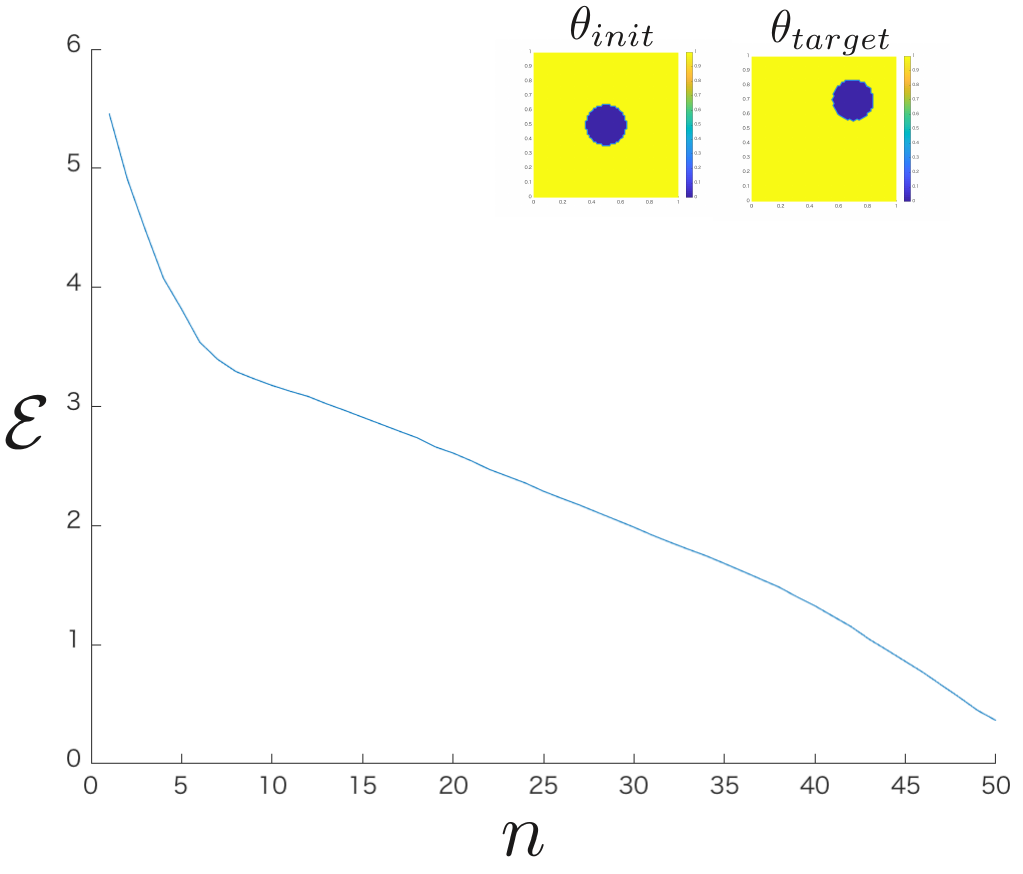}
  \caption{Evolution of the cost functional $\mathcal{E}$ under the gradient flow. Initial and target conditions are also displayed.}
  \label{fig:graph4}
\end{figure}

\begin{figure}[H]
  \centering
  \includegraphics[scale=0.5,page=6]{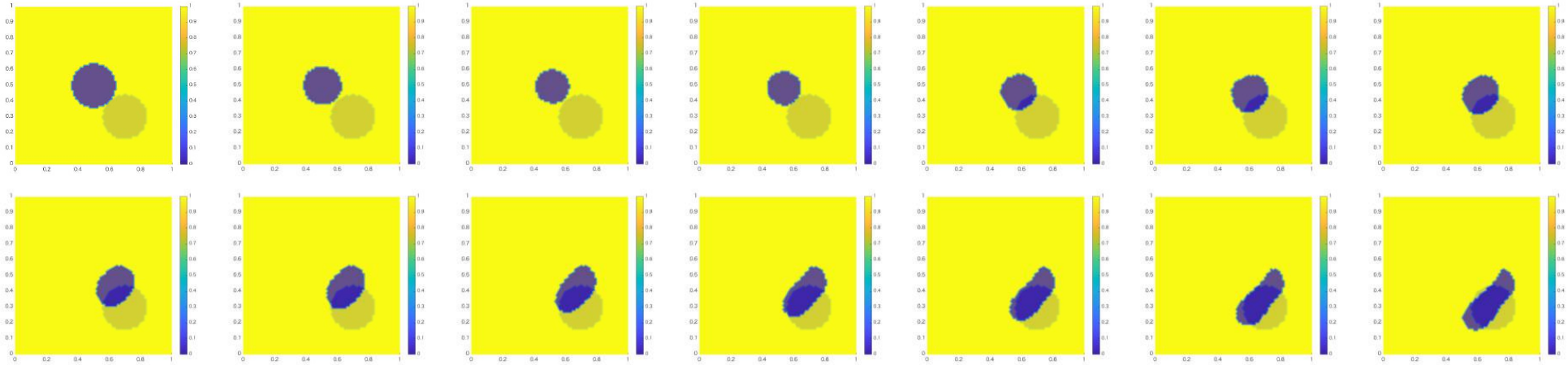}
  \caption{Numerical solution of the inverse problem. Time is from top to bottom, left to right. The dark region is the approximate solution and $\theta_{target}$ is shown as the light disk.}
  \label{fig:optim8}
\end{figure}

\begin{figure}[H]
  \centering
  \includegraphics[scale=0.4,page=5]{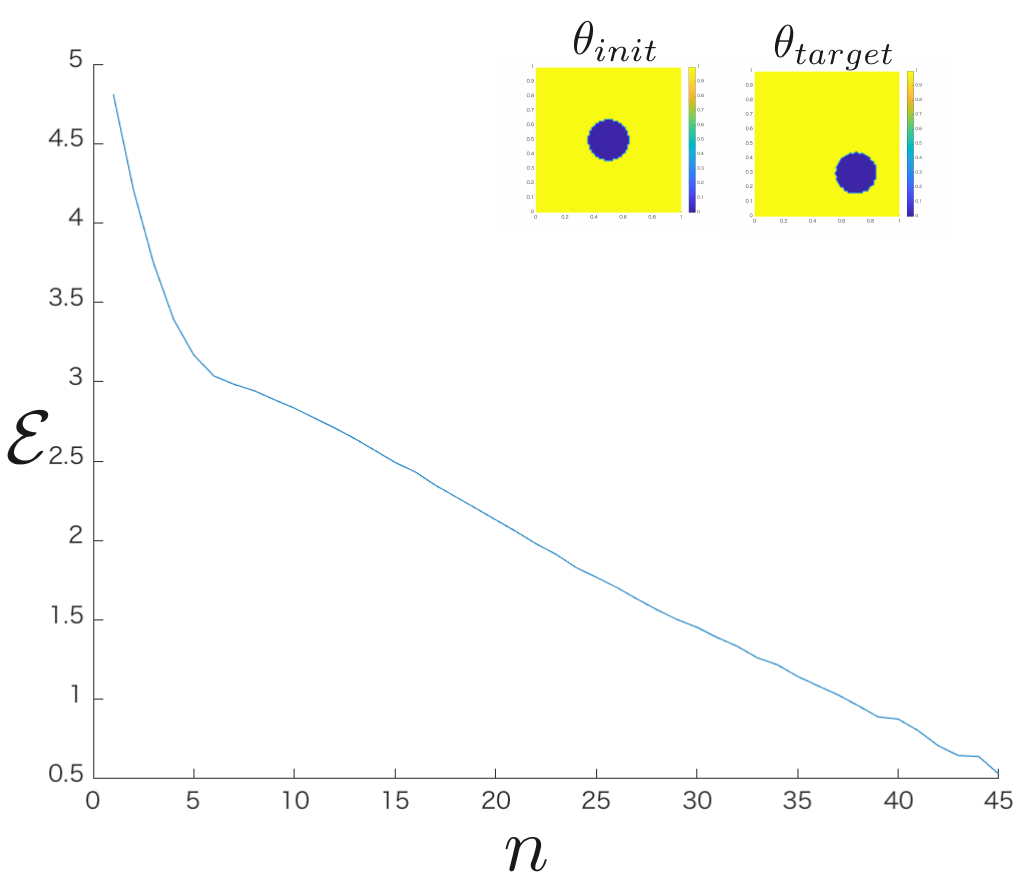}
  \caption{Evolution of the cost functional $\mathcal{E}$ under the gradient flow. Initial and target conditions are also displayed.}
  \label{fig:graph5}
\end{figure}

We also examine the gradient flow of the cost functional when the target is placed in the upper right corner of the computational domain (as in figure \ref{fig:g_confirm} b). The evolution is shown in figure \ref{fig:optim7}. The level set function converges to the approximate location of $\theta_{target}$ while decreasing the value of the cost functional (see figure \ref{fig:graph4}). Although the shape of the limit function is not as precise as those obtained in our the previous computations, we remark that its overall shape and location is correct. We also recall that these results are obtained through the use of the boundary data $\delta$, and nothing more. Boundary conditions could be effecting the computational result, but we also expect that mesh and parameter dependence may be playing a role as well. Similar results are obtained when the target is place in the bottom right corner of the computational domain (see figures \ref{fig:optim8} and \ref{fig:graph5}).

\begin{figure}[H]
  \centering
  \includegraphics[scale=0.5,page=5]{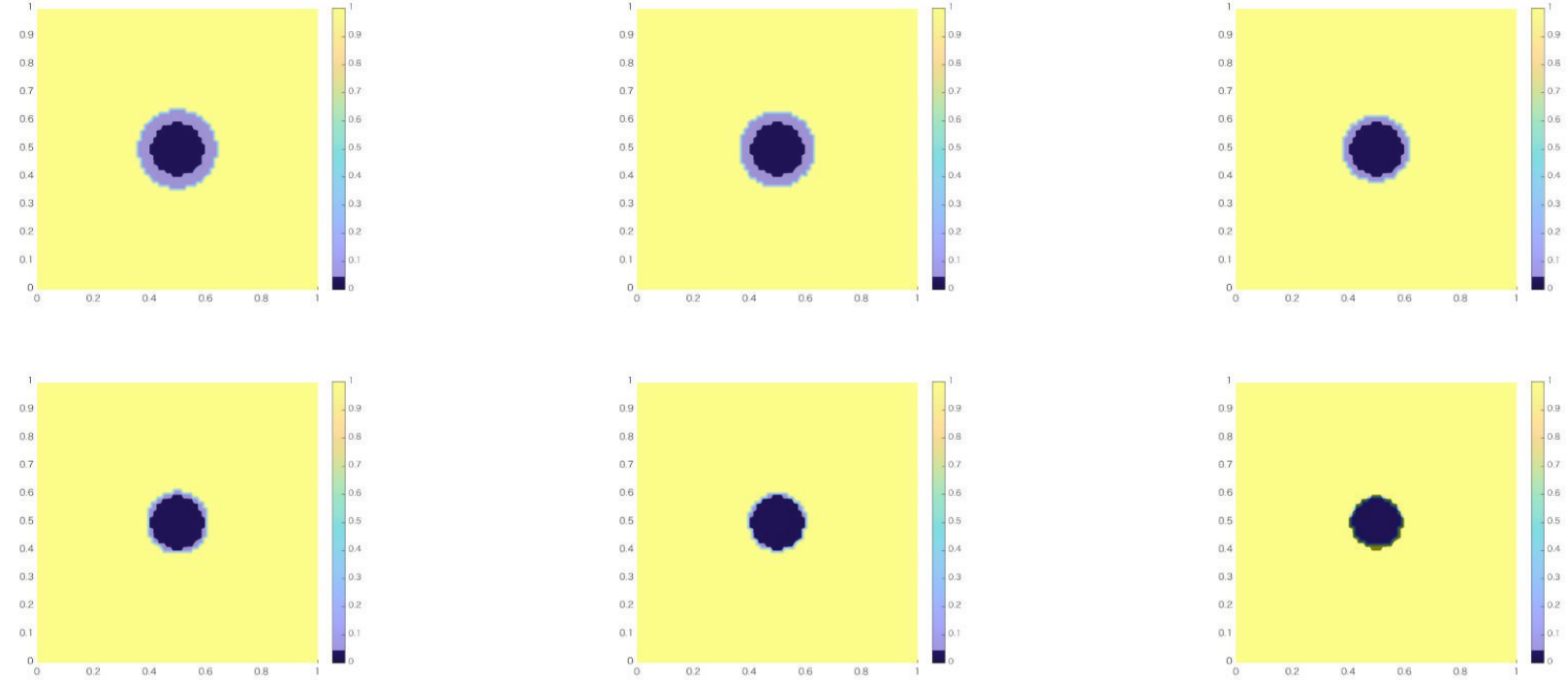}
  \caption{Numerical solution of the inverse problem. Time is from top to bottom, left to right. The dark region denotes the target.}
  \label{fig:optim9}
\end{figure}

\begin{figure}[H]
  \centering
  \includegraphics[scale=0.4,page=6]{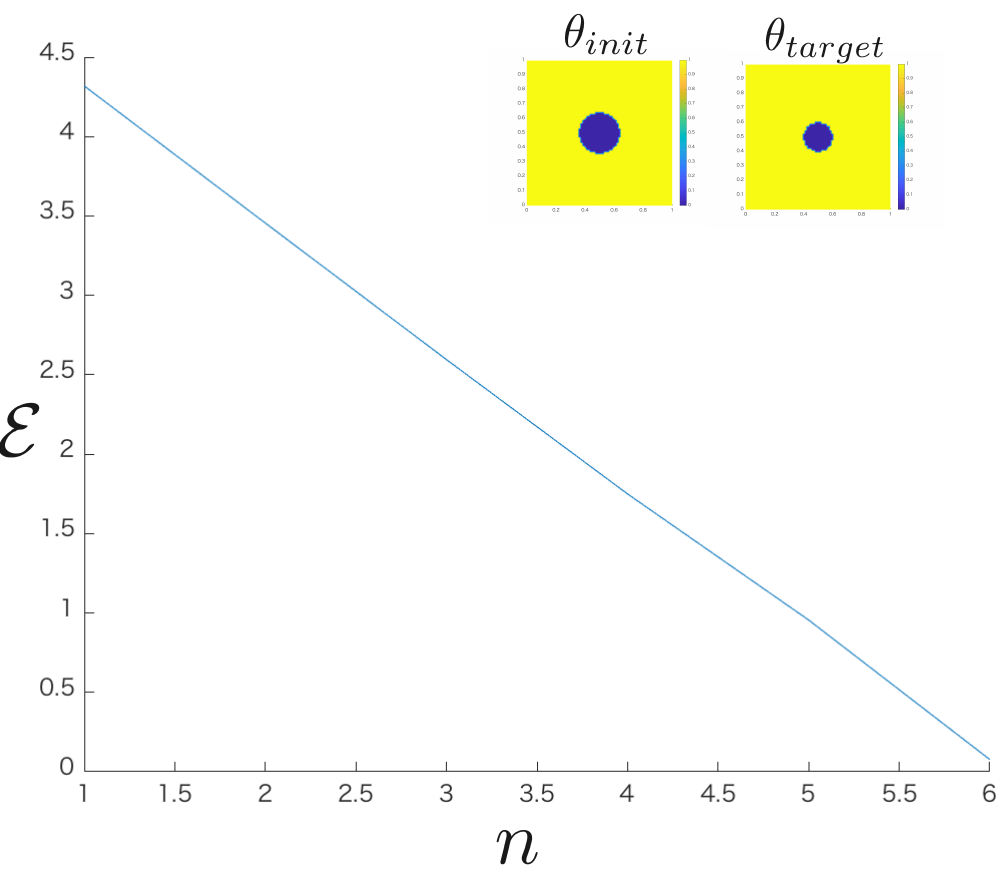}
  \caption{Evolution of the cost functional $\mathcal{E}$ under the gradient flow. The initial condition is shown overlapping the target (shown in dark).}
  \label{fig:graph6}
\end{figure}

We also considered the case where $\theta_{init}$ overlaps $\theta_{target}$ with a larger radius. The numerical results are shown at six instances in figure \ref{fig:optim9}. Here, we observe the level set function to shrink inward while converging to its target. The evolution of the cost functional's values under the gradient flow are displayed in figure \ref{fig:graph6} which shows a linear decrease.

\noindent

\begin{figure}[H]
  \centering
  \includegraphics[scale=0.5]{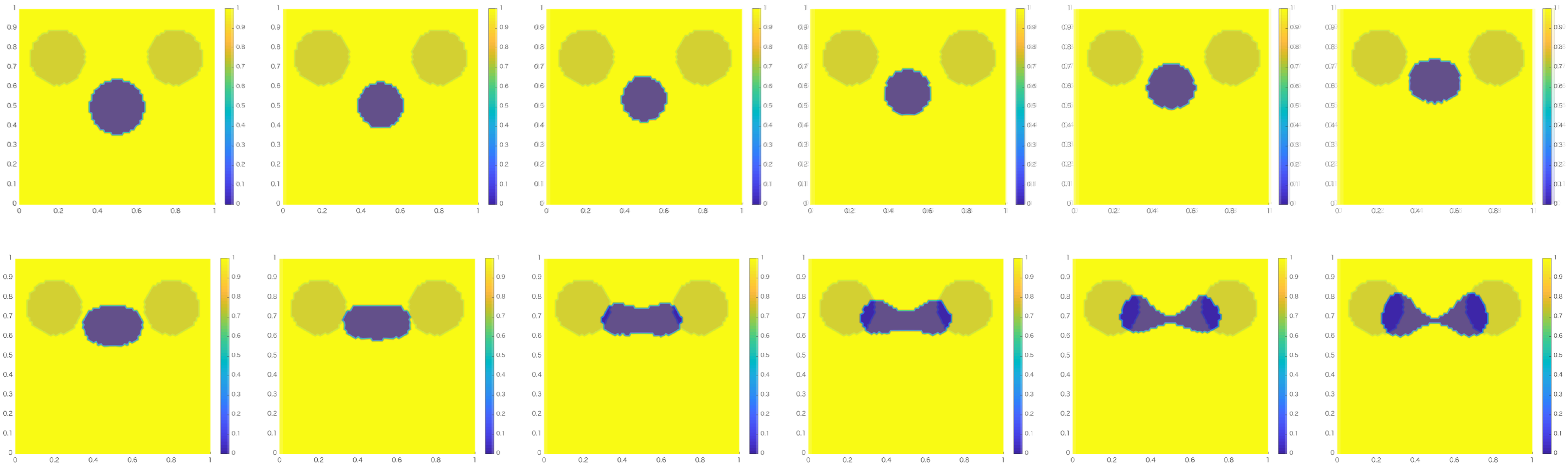}
  \caption{Evolution of the level set function with a disconnected target (shown as a shadow). Time is from top to bottom, left to right.}
  \label{fig:future1}
\end{figure}

Our final numerical experiment sets the target as two disconnected disks with radii equal to that of $\theta_{init}$ (see figure \ref{fig:future1}). The numerical results show the level set function spread towards the target and finally splits into two disconnected regions. Although the cost functional's value decreases with each iteration of our approximation method, we note that the limit function is visibly different from the target. Nevertheless, the limit function visually does represent an approximation of the target and this fact suggests the possible existence of local minimizers of the cost functional near the target.

The computations in this section clearly show the ability of our method to approximate solutions of the inverse problem stated in section \ref{inv}. Overall, these results show that the cost functional's gradient indeed contains information about the interior composition of the domain.

\begin{figure}[H]
  \centering
  \includegraphics[scale=0.5]{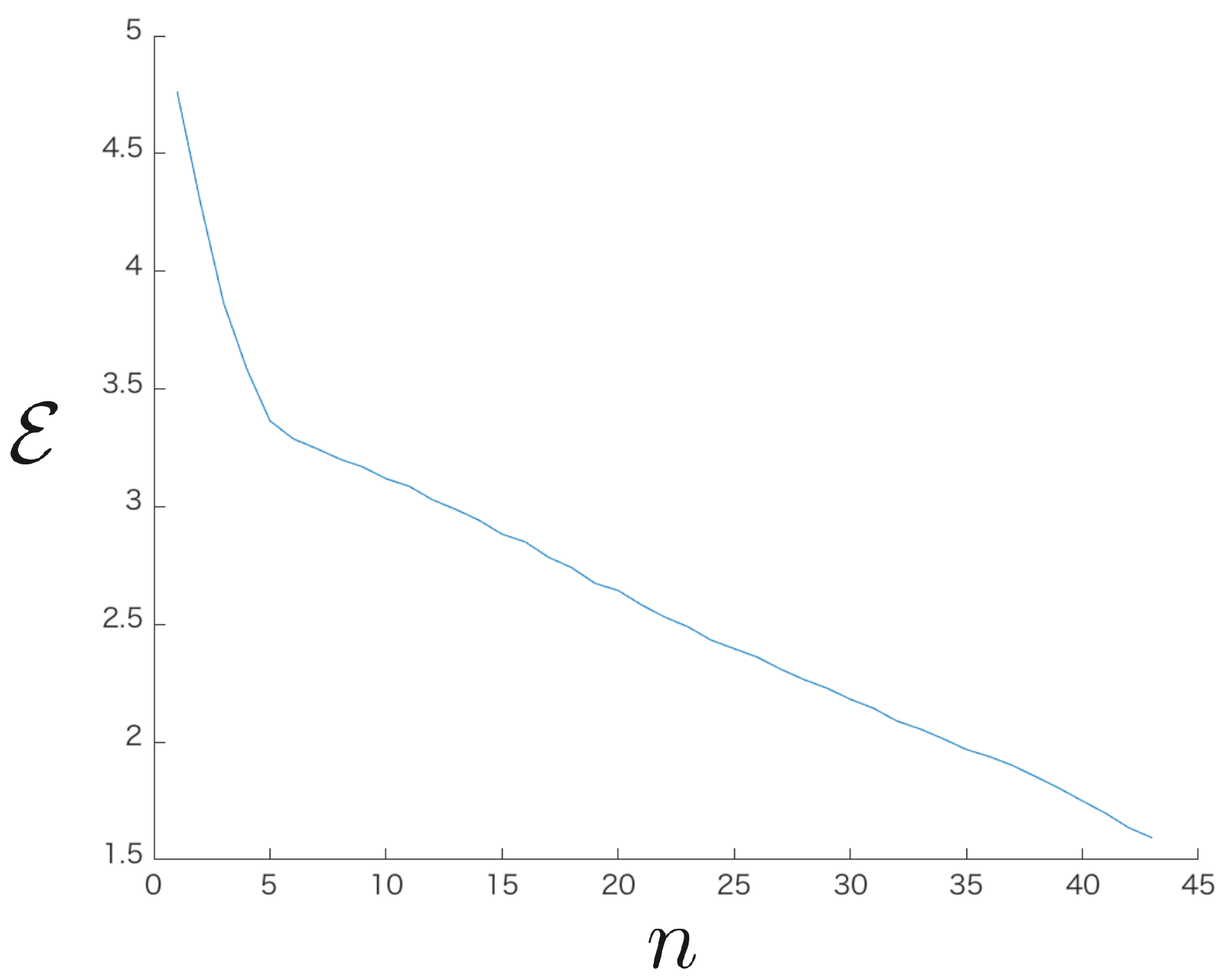}
  \caption{Evolution of the cost functional under a gradient flow towards the disconnected target.}/
  \label{fig:future2}
\end{figure}

\section{Conclusion and future work}\label{conclusion}

A mathematical technique for determining the internal structure of composite elastic bodies was developed. We showed how a single component of the displacement field along a portion of the boundary of the elastic body can be used to gain information about the interior composition. This was done by using the method of Lagrange multipliers to derive the gradient of a cost functional whose minimization corresponds to our target inverse problem.

Numerical solvers were also developed for computing solutions to model and adjoint equations. These were used to construct numerical realizations of the cost functional's gradient. In addition, we designed a level set method which incorporates the cost functional's gradient to evolve the level set function as the gradient descent of the cost functional. This allowed us to examine the numerical solution of the inverse problem and our results showed that our framework is able to recover the interior composition of composite elastic materials, using only a single component of the displacement field's boundary data. 

We would like to continue improving our framework for use on problems using real experimental data, such as the surface acoustic wave data mentioned in the introduction (see also \cite{Ohtsuka}). Here we are particularly interested in whether or not our method can be used to image the interior of composite elastic bodies which contain sub-micron scale inclusions. Relatedly, we aim to develop numerical methods for treating the three dimensional phenomena as well. Since the numerical methods developed in this study tend to be computationally heavy, we would like to improve their efficiency. Finally, we are also planning to investigate inverse problems involving a larger number of composite materials.

\section{Acknowledgements}
The work of E. Ginder was supported by JST Presto ``Multiphase shape optimization in phononic crystal design." The author would also like to acknowledge support from JSPS KAKENHI Grant Number  17K14229 and JSPS KAKENHI Grant Number 15KT0099.

\end{document}